\documentclass{birkjour_t2}
\usepackage{amssymb,amsmath,ma 
   thrsfs,amsfonts,amsthm,cases}

\usepackage{enumerate}
\usepackage{xcolor}
\usepackage[numbers]{natbib}
\usepackage{geometry}
\usepackage{setspace}
\usepackage{appendix}
\usepackage[colorlinks,
            linkcolor=teal,
            anchorcolor=blue,
            citecolor=purple]{hyperref}

\renewcommand{\labelenumi}{(\roman{enumi})}

\begin{document}

\newtheorem{theorem}{Theorem}[section]

\newtheorem{lemma}[theorem]{Lemma}
\newtheorem{prop}[theorem]{Proposition}
\newtheorem{cor}[theorem]{corollary}
\newtheorem{remark}[theorem]{Remark}
\newtheorem*{con}{Conjucture}
\newtheorem{assumption}{Assumption}
\newtheorem*{questionA}{Question}
\newtheorem{definition}{Definition}[section]

\def \elam {\eta^\nu}
\def \ulam {u^\nu}
\def \plam {\psi^\nu}
\def \ru {\mathcal{R}_{u}}
\def \ret {\mathcal{R}_{\eta}'}
\def \rqut {\mathcal{R}_{\qu}'}

\def \pulam {\P u^{\nu}}
\def \qulam {\q u^{\nu}}
\def \ett {e^{\frac{2\tc}{ \nu}t}} 
\def \ets {e^{\frac{2\tc}{ \nu}s}} 
\def \et {e^{\frac{\tc}{ \nu}t}} 
\def \es {e^{\frac{\tc}{ \nu}s}}

\def \ne {\nu^2}
\def \ney {\nu^{\frac{1}{2}}}
\def \nfey {\nu^{-\frac{1}{2}}}
\def \ny {\nu^{-1}}
\def \nfe {\nu^{-2}}
\def \nes {\nu^{-\frac{3}{2}}}

\def \g {\nabla}
\def \gyi {{\nabla^1}}

\def \curl {\g\times}
\def\p{\partial}
\def \A {\cal A}
\def \le {\lesssim}

\def \bl {{\big\{}}
\def \brr {{\big\}}}
\def \bbl {{\Big\{}}
\def \bbr {{\Big\}}}

\def \div {\text{div}\,}
\def \det {\text{det}}
\def \divr {\text{div}_R}
\def \d {\,\mathrm{d}}
\def \dx {\,\mathrm{d}x}
\def \dr {\,\mathrm{d}R}
\def \dri {\,\frac{\mathrm{d}R}{\psi_\infty}}
\def \dt {\frac{\d }{\d t}}
\def \drx {\dri\dx}
\def \er {\frac{\eta}{\rho}}
\def \geta {{\g\eta }}
\def \etm {{\eta^m}}
\def \qum {{\qu^m}}
\def \pum {{\pu^m}}
\def \um {{u^m}}

\def \pru {{\P \cal M}}
\def \eu {{(\eta u)}}

\def \emm {\eta^{m-1}}
\def \qumm {\qu^{m-1}}
\def \pumm {\pu^{m-1}}
\def \etmm {{\eta^{m-1}}}
\def \getam {{\g\eta^{m-1} }}
\def \umm {u^{m-1}}
\def \wmm {w^{m-1}}

\def \gr {\g_R}
\def \div {\text{div}\,}
\def \divr {\text{div}_R}
\def \gm {\g^m}
\def \gmp {\g^{m+1}}
\def \gmm {\g^{m-1}}
\def \mm {^{m-1}}

\def \eb {\cal E_B}
\def \ep {\cal E_P}
\def \ei {\cal E_I}
\def \ee {\cal E_\eta}

\def \emr {\eta^{m-r}}
\def \qumr {\qu^{m-r}}
\def \pumr {\pu^{m-r}} 
\def \umr {u^{m-r}}

\def \r {\mathbb{R}}
\def \T {\mathbb{T}}
\def \D {\mathbb{D}}
\def \P {\mathbb{P}}
\def \q {\mathbb {Q}}
\def \h {\mathcal}
\def \hlt {{\h L^2}}
\def \hd {{{\dot{\mathcal H}}^1}}
\def \hdt {{\dot{\mathcal H}}^2}
\def \pin {\psi_\infty}
\def \vp {\Psi}
\def \buk {U^\kappa}
\def \uk {u^\kappa}
\def\bpk{\Psi^\kappa}
\def\pk{\psi^\kappa}

\def \tbuk {\lt {U^\kappa}}
\def \tuk {\lt {u^\kappa}}
\def \abuk {\lb {U^\kappa}}
\def \auk {\lb {u^\kappa}}

\def \tc { \tilde c }
\def \bc {\bar C}

\def \cu {\mathbf}
\def \cal{\mathcal}

\def \t {\tilde}
\def \b {\bar}

\def \pu {\mathbb P u}
\def \qu {\mathbb Q u}

\def \rml { \rho^{\mu,\lambda}}
\def \eml { \eta^{\mu,\lambda}}
\def \uml { u^{\mu,\lambda}}
\def \quml { \qu ^{\mu,\lambda}}
\def \puml { \pu ^{\mu,\lambda}}

\def \el {{\eta^L}}
\def \ul {{\qu^L}}

\def \ler {{\h L^2}}
\def \lerer {{\h L^2}^2}
\def \hyi {{\dot {\h H}^1}}
\def \hyier {{\dot {\h H}^1}^2}
\def \sler {s,{\h L^2}}
\def \slerer {{s,{\h L^2}}^2}
\def \mler {{m,{\h L^2}}}
\def \mlerer {{m,{\h L^2}}^2}
\def \mmler {m-1,{\h L^2}}
\def \mmlerer {{m-1,{\h L^2}}^2}
\def \shyi {{s,\hyi}}
\def \shyier {{s,\hyi}^2}
\def \mhyi {{m,\hyi}}
\def \mhyier {{m,\hyi}^2}
\def \mmhyi {{m-1,\hyi}}
\def \mmhyier {{m-1,\hyi}^2}

\def \hm {{H^m}}
\def \hmm {{H^{m-1}}}
\def \hmm {{H^{m-1}}}

\def \lyi {{L^1}}
\def \lyier {{L^1}^2}
\def \m {\cal M(t)}
\def \mer {\cal M^2(t)}

\def \ix {\int_{\T^3}}
\def \it {\int _0^t}

\def \supt {\sup\limits_{0\leq s\leq t}}

\def \intr {\int_B}
\def \ixr {\iint_{\T^3\!\times \!B}\!}

\def \ph {\pin^{\frac{1}{2}}}
\def \phf {\pin^{-\frac{1}{2}}}

\def \rqu {\mathcal{R}_{\qu}}
\def \rpu {\mathcal{R}_{\pu}}
\def \rqum {\mathcal{R}^m_{\qu}}
\def \rem {\mathcal{R}^m_{\eta}}
\def \rqumm {\mathcal{R}^{m-1}_{\qu}}
\def \remm {\mathcal{R}^{m-1}_{\eta}}
\def \rwmm  {\mathcal{R}^w_{\qu}}

\def \holder {H$\ddot{\text{o}}$lder's }
\def \poin {Poincar$\Acute{\text{e}}$ }
\def \m {{\mathfrak m}}

\bibliographystyle{abbrvnat}

\title[Incompressible limit of FENE dumbbell]{Incompressible limit for the 3D compressible FENE dumbbell model}


\author[J. Gao]{Jincheng Gao}
\address{
School of Mathematics
Sun Yat-sen University
Guangzhou  510275
People’s Republic of China\br
email: gaojch5@mail.sysu.edu.cn}

\author[J. Wu]{Jiahong Wu}
\address{
Department of Mathematics,
University of Notre Dame,
Notre Dame, IN 46556
USA\br 
email:jwu29@nd.edu}

\author[Z. Yao]{Zheng-an Yao}
\address{
School of Mathematics
Sun Yat-sen University
Guangzhou  510275
People’s Republic of China \br 
email: mcsyao@mail.sysu.edu.cn}

\author[R. Yu]{Ruijia Yu$^*$}
\address{ 
School of Mathematics
Sun Yat-sen University
Guangzhou  510275
People’s Republic of China \br
email: yurj5@mail2.sysu.edu.cn}
\begin{abstract}
The FENE dumbbell model couples the Navier–Stokes equations governing fluid velocity with a Fokker–Planck equation describing the evolution of polymer distribution within the fluid.
In this work, we study the global-in-time incompressible limit of the compressible FENE dumbbell model on the three-dimensional torus $\T^3$, where the incompressible limit is driven by large volume viscosity. To establish this limit, we develop time-weighted a priori estimates that yield decay rates for strong solutions. A key challenge arises from the fact that increasing the volume viscosity suppresses the decay of high-frequency components, thereby weakening the dissipation of the density and complicating the derivation of uniform-in-time decay estimates. To overcome this difficulty, we introduce a novel momentum-based estimate and show that the incompressible component of the momentum decays faster in time than the velocity itself. Exploiting this enhanced decay, we successfully close the a priori estimates and establish a time-decreasing convergence rate toward the incompressible limit.
\end{abstract}


\maketitle
\section{Introduction}
Partial differential equation (PDE) systems modeling the interaction between fluids and polymers have received growing attention in recent years due to their wide-ranging applications in physics, chemistry, and biology \cite{bird1, masmoudi-cpam}.
Among these models, the finite extensible nonlinear elastic (FENE) dumbbell model stands out as a particularly significant example. In this framework, a polymer is idealized as an “elastic dumbbell,” consisting of two beads connected by a nonlinear spring, with the configuration described by a vector $R$ \cite{bird1, bird2, doi, ottinger}. At the fluid level, the FENE dumbbell model couples the Navier–Stokes equations, which govern the fluid velocity, with a Fokker–Planck equation that describes the evolution of the polymer distribution within the fluid medium. The micro–macro compressible FENE dumbbell model takes the following form \cite{bird1, doi}:
\begin{equation}\label{fene-orgin}
    \begin{cases}
        \p_t \rho + \div (\rho u)=0,\\
        \p_t (\rho u) + \div (\rho u\otimes u) -\mu \Delta u -( \mu + \lambda)\g\div u +\g P = \div \tau,\\
          \p_t \Psi + u\cdot \g \Psi
          = \divr\left[
          -\sigma (u)\cdot R\Psi+\beta\gr \Psi +\gr\, \mathcal U \Psi
          \right ],\\
          \tau_{j,k}=\intr (R_j \p_{R_k} \mathcal U)\Psi(x,R,t)\dr,\\
          \rho (x,0)=\rho_0,\quad
          u(x,0)=u_0,\quad \Psi(x,R,0)=\Psi_0,\\
          (\gr \Psi +\gr\, \mathcal U \Psi)\cdot n = 0 \quad \text{on} \ \p B(0,R_0).
    \end{cases}
\end{equation}
In system \eqref{fene-orgin}, $\rho(x,t)$ denotes the fluid density, $u(x,t)$ is the velocity of the polymeric liquid, $P(x,t)$ is the pressure, and $\Psi(x,R,t)$ represents the distribution function of the internal polymer configuration. Here, the spatial variable $x$ lies in the three-dimensional torus $\T^3 := [-\pi, \pi]^3$, and the polymer elongation vector $R$ is constrained within the ball $B(0, R_0) \subset \mathbb{R}^3$, reflecting the finite extensibility of the polymers. To simplify the presentation, we assume
 $$
 P(\rho)=\frac{\rho^3}{3},
 $$
noting that the more general case $P(\rho) = a\rho^\gamma$ with $a > 0$ and $\gamma > 1$ can be addressed using Taylor’s formula near the constant background. The potential $\h U$ governing the polymer elasticity is given by
 $$\h U(R)=-k \log \left( 1-\frac{|R|^2}{R_0^2} \right)$$
  for some $k>0$. The magnitude of $k$ is an important parameter as it portrays the strength of singularity of the equilibrium state $\pin$ on the boundary $\p B(0,R_0)$.
In addition, $\mu$ and $\lambda$ are the shear viscosity and volume viscosity of the fluid, respectively. $\beta$ relates to the Boltzmann constant and temperature. 
  In general, $\sigma(u)=\g u$. For the co-rotation case, $\sigma(u)=\frac{\g u-(\g u)^{\text{T}}}{2}$. For simplicity, we will set $\beta=1$ and $R_0=1$. 

  As in \cite{masmoudi-cpam}, to ensure the conservation of $\Psi$, we add an additional boundary condition, namely
  $$
  (-\g u\cdot R\Psi+\beta\gr \Psi +\gr\, \mathcal U \Psi)\cdot n=0 \quad \text{ on }\p B(0,R_0).
  $$
  This boundary condition implies that $\Psi=0$ on $\p B(0, R_0)$, and if $\intr \Psi_0( x, R ) \dr=1$, then for all $t$, we have $\intr \Psi( x,R,t ) \dr=1$.
  
 
In compressible flows, the volume (or bulk) viscosity can significantly exceed the shear viscosity. This is observed in gases such as carbon dioxide, methane, and nitrous oxide \cite{sharma}. The effect is especially pronounced during compression or expansion, where the relaxation time required to restore thermodynamic equilibrium is long, leading to a substantial increase in bulk viscosity \cite{landau}. 

High-volume-viscosity fluids arise in a variety of practical applications, including power systems driven by non-fossil fuel heat sources, aerodynamic testing in wind tunnels, and numerous processes in pharmaceutical manufacturing \cite{cramer2012numerical}.

Motivated by these physical and practical considerations, we study the stability and long-time behavior of system \eqref{fene-orgin} in the incompressible limit as $\lambda \to \infty$.
Our analysis crucially relies on the formal observation that, as  $\lambda\rightarrow \infty$, the compressible system \eqref{fene-orgin} converges to the incompressible FENE dumbbell model:
 \begin{equation*}
     \begin{cases}
        \p_t v + v\cdot\g v - \mu \Delta v + \g\Pi=\div \t \tau,\quad \div v=0,\\
        \p_t \varPhi + v\cdot\g \varPhi=\divr\left[
          -\g v\cdot R\varPhi+\beta\gr \varPhi +\gr\, \mathcal U \varPhi
          \right ],
     \end{cases}
 \end{equation*}
 where $\Pi$ is the pressure.
 
The well-posedness of various forms of the FENE dumbbell model has been extensively studied. For the incompressible FENE dumbbell system,
\citet{renardy} established local well-posedness in Sobolev spaces for potentials of the form $\widehat{U}(R) = (1 - |R|^2)^{1 - \sigma}$ with $\sigma > 1$.
Subsequently, \citet{jourdain} proved the local existence of solutions to the corresponding stochastic differential equation with potential $\widehat{U}(R) = -k \log(1 - |R|^2)$, assuming $b = 2k > 6$, in the setting of Couette flow. \citet{zhang-arma} further extended the analysis to three dimensions, establishing local well-posedness in weighted Sobolev spaces when $b = 2k > 76$. A major breakthrough was achieved by \citet{masmoudi-cpam}, who introduced novel Hardy-type inequalities to handle the singular term $\operatorname{div} \tau$, and proved both local and global well-posedness near equilibrium for the FENE model under the mild assumption $b = 2k > 0$. In addition, \citet{masmoudi-invension} proved the global existence of weak solutions in $L^2$, under suitable entropy conditions. More recently, \citet{lin-cpam} established the global existence of solutions near equilibrium under certain structural constraints on the potential.

For the compressible FENE dumbbell models, 
\citet{barrett-m3} established the existence of global-in-time weak solutions to the compressible FENE dumbbell models with center-of-mass diffusion in two or three-dimensional bounded domains with pressure $P(\rho)=C\rho^\gamma$ and $\gamma > \frac{3}{2}$. 
Their analysis required the inclusion of a quadratic interaction term in the elastic extra-stress tensor, in addition to the classical Kramers expression, to achieve the compactness necessary for the proof.
In contrast, for the two-dimensional case with $\gamma>1$, \citet{barrett-jde} established the existence of global-in-time weak solutions without the need for the quadratic interaction term.
 \citet{breit} considered the compressible FENE dumbbell model with classical Warner potential, and obtained the unique local existence of the weak solution in $\T^2$ or $\T^3$ with or without the center-of-mass diffusion.
Very recently, \citet{luo2024} proved the global existence of the compressible FENE dumbbell model in $\r^d$ with $d\geq 2$ under the Besov framework.

Although extensive research has been conducted on both compressible and incompressible FENE dumbbell models, to the best of our knowledge, there is no existing work that establishes a direct connection between these two systems (for dilute models, see \cite{suli2020}). 
To address this gap, we investigate the relationship between the compressible and incompressible FENE dumbbell models via the incompressible limit process.

   The incompressible limit driven by large volume viscosity has been the subject of extensive study in recent years.
   \citet{danchin2017} established the global existence and incompressible limit of compressible Naiver-Stokes equation in $\r^2$, considering large initial velocity in critical Besov space. They also extended their analysis to $\r^3$,  assuming the existence of solutions.
  Later, \citet{cui-3d} employed the vector field method to investigate the incompressible limit of compressible viscoelastic systems when the shear viscosity converges to zero.
  Subsequently, \citet{cui-2d} addressed the challenges posed by weaker decay rates and extended their previous results to $\r^2$. 
  For other related results on this type of incompressible limit, see \cite{chen2019,danchin2018-t, hu2023incompressible} and the references therein.
It is worth noting that, apart from large volume viscosity, another mechanism for attaining the incompressible limit is the low Mach number approach. 
Low Mach number induced incompressible limit of compressible Navier-Stokes equations is widely studied in whole space or periodic domains \cite{alazard,desjardins,feireisl2007,hu2009,huang2017, levermore} and bounded domains \cite{jiang2011,ju,klainerman, masmoudi2022,ou2022}.
Also, for low Mach limit of compressible MHD systems, please refer to \cite{hu2009,jiang2014low,klainerman, majda2012compressible}.

To close the uniform estimates for the incompressible limit, the $L^1$-time integrability of the incompressible part of velocity is needed (see Section \ref{sub-difficulties} for details). In \cite{danchin2017}, this is achieved by employing the Besov space framework, where $L^1$-time integrability follows naturally. 
 However, within the Sobolev framework, additional analysis of the decay properties is required to obtain the desired $L^1$-time integrability.
Moreover, substantial research has been devoted to exploring the decay of the FENE dumbbell models.
For incompressible FENE dumbbell models,
 \citet{schonbek} studied the $L^2$ decay of the velocity for the co-rotation FENE dumbbell model, and obtained the decay rate $(1+t)^{ -\frac{N}{4}+\frac{1}{2} }$.
 Later, Luo and Yin \cite{luo-arma,luo-advance} improved the $L^2$ decay results developed in \cite{schonbek}  by Fourier splitting methods, and obtained the optimal decay rates for the co-rotation case.
While for compressible viscoelastic dumbbell models, \citet{luo2024} proved the optimal decay rates for $d\geq 3$ by the Littlewood-Paley decomposition theory and the Fourier splitting method.
 There are also many important results on the decay of other viscoelastic models, for instance, see \cite{tan-jde,chen-siam,hu-sima,huang2022,wang2022sharp} and the references therein.

\subsection{Statement of results}
In this paper, we first investigate the global stability and large-time behavior of (\ref{fene-orgin}), then consider the convergence rate of (\ref{fene-orgin}) under the limit process $\lambda\rightarrow \infty$. In particular, we focus on the dynamics near the equilibrium $\rho=\bar\rho=1$, $u=0$, $\Psi=\pin$, where
  $$
  \bar \rho=\ix\rho_0 \d x,\qquad
  \pin(R)=\frac{e^{-\h U(R)}}{\intr e^{-\h U(R)}
\dr }=\frac{(1-|R|^2)^k}{\intr (1-|R|^2)^k \dr}.
  $$
For this purpose, we denote $\nu:=2\mu + \lambda$, $\elam := \rho-\bar\rho$, $\rho^\nu:=\rho$, and $\plam :=\Psi-\pin$.
 Since $\g_x \pin=0$, we have 
 $$
 \div \tau=\div \intr (R\otimes\gr \h U)\Psi \dr = \div \intr (R\otimes\gr \h U)\plam \dr.
 $$
 Hence, we may assume that
 $$
 \tau= \intr (R\otimes\gr \h U)\plam \dr.
 $$
  
With the above notations, the equation governing the perturbation $(\elam, \ulam,\plam)$ of the $3$D compressible FENE dumbbell model reads as follows:
      \begin{equation}\label{fene}
      \begin{cases}
          \p_t \elam + \ulam\cdot\g \elam +\div \ulam= -\elam \div \ulam,\\
          \p_t \ulam + \ulam\cdot\g \ulam -\mu\Delta \ulam -( \mu + \lambda )\g\div \ulam + \g \elam =\div \tau^\nu +  \ru,\\
          \p_t\plam +\ulam\cdot\g\plam =
          \divr(-\g\ulam\cdot R (\plam+\pin)) + \divr(\pin\gr\frac{\plam}{\pin}),
      \end{cases}
  \end{equation} 
  where
  $$
  \ru = -\mu\frac{\elam}{\rho^\nu}\Delta \ulam -( \mu + \lambda )\frac{\elam}{\rho^\nu}\g\div \ulam - \elam\g\elam
  -\frac{\elam}{\rho^\nu}\div \tau^\nu.
  $$

Since momentum plays a pivotal role in our analysis, we denote it as $\mathcal{M}^\nu := \rho^\nu u^\nu$.
Apart from the regularity assumption, we also assume that 
\begin{equation}\label{0meanvalue}
    \ix \mathcal{M}^\nu (0) \dx=0,
\end{equation}
Since the Leray Projector $\P$ does not affect the zero mode in frequency space, for all $t>0$, we have
\begin{equation}\label{mean value 0}
    \ix \elam(t) \dx=0,\quad
    \ix \mathcal{M}^\nu (t) \dx
    =\ix \pru^\nu (t) \dx=0.
\end{equation}

Our first main result considers the global existence and decay property of (\ref{fene}).
 \begin{theorem}\label{thm-exist}
    Let $m\geq 3$ be an integer. Assume that the initial data satisfies
    $$
    (\eta_0, u_0,\psi_0)\in H^m(\T^3)\times\{H^m(\T^3)\}^3\times H^m(\T^3,\ler).
    $$ 
    Moreover, assume that (\ref{0meanvalue}) holds. Then there exist a small constant $\varepsilon>0$ and a large constant $\nu_0>0$ such that if $\nu\geq \nu_0$ and
    $$
     ||\ulam_0||_m + \ney||\qulam_0||_m
     + \ney||\elam_0||_m + ||\plam_0||_{m,\h L^2}
     \leq \varepsilon,
    $$
    system (\ref{fene}) admits a unique global classical solution 
    $$(\elam, \ulam,\plam)\in H^m(\T^3)\times\{H^m(\T^3)\}^3\times H^m(\T^3,\ler).
    $$
    Moreover, for any $t>0$, the following estimates hold: 
    $$
    \begin{aligned}
    &\supt \ets || \ulam(s)||_{m}^2
    + \int_0^t\ets \bbl  ||\elam(s)||_{m}^2 + \mu || \ulam (s)||_{m+1}^2 \bbr \d s\leq C \varepsilon^2, \\
    &\supt\bl  \nu\ets ||\qulam(s)||_{m}^2
    +  \nu\ets  ||\eta(s)||_{m}^2 \brr
    + \nu^2\int_0^t \ets ||\qulam(s)||_{m+1}^2\d s\leq C \varepsilon^2,   \\
    &\supt \bl  a(s) 
    ||\pru^\nu(s)||_{m-1}^2 
    +  a(s) ||\plam(s)||_{m-1,\ler}^2
    \brr
    \leq C \varepsilon^2,\\
    &\int_0^t a(s)
    \bbl ||\pru^\nu(s)||_{m}^2
    + ||\plam(s)||_{m-1,\hd}^2\bbr \d s \leq C \varepsilon^2,
    \end{aligned}
    $$
    where     
    \begin{equation}\label{at}
        a(t)=
        \begin{cases}
            (1+\delta t )^2, \qquad &0< t \leq \nu,\\
            C_\delta \nu^2\ett,\qquad &t> \nu,
        \end{cases}
    \end{equation}
    and $\tc$, $\delta$ are small positive constants, and $C_\delta$ is a positive constant related to $\delta$.
\end{theorem}   

\begin{remark}
    We can immediately deduce from the above theorem that $(\elam, \ulam,\plam)$ satisfies the following decay properties :
    \begin{equation*}
        \begin{aligned}
     &||\ulam(t)||_m\leq C \varepsilon e^{-\frac{\tc}{\nu }t},\qquad
     ||\qu^\nu(t)||_m
     \leq C \varepsilon\nfey e^{-\frac{\tc}{\nu }t},\\       
    &||\P\ulam(t)||_m
    \leq C \varepsilon
    \max\{
    \nfey e^{-\frac{\tc}{\nu }t},
    a(t)^{-\frac{1}{2}}\},\qquad
    ||\pru^\nu (t)||_{m-1}
    \leq C \varepsilon a^{-\frac{1}{2}}(t),\\
    &||\elam(t)||_m
    \leq C \varepsilon\nfey e^{-\frac{\tc}{\nu }t},\qquad
     ||\plam(t)||_{m-1,\ler }
    \leq C \varepsilon a^{-\frac{1}{2}}(t).
        \end{aligned}
    \end{equation*}
A key and intriguing observation in this study is that, the incompressible part of the momentum $\mathbb{P}(\rho u)$ exhibits faster time decay rate compared to $\mathbb{P} u$ (please refer to Section \ref{sub-difficulties} for more details).
 This observation not only yields an improved decay rate but also offers valuable insights and a novel perspective for future investigations into fluid models with large volume viscosity.
\end{remark}

\begin{remark}
     Although the exponential decay rate involves the parameter $\frac{1}{\nu}$ and is therefore not uniform as $\nu\rightarrow \infty$, removing this parameter is not expected.
    If the effect of the elastic dumbbell is neglected, namely $\psi=0$, then (\ref{fene}) becomes the classical compressible Navier-Stokes equation. Classical decay results (see e.g. \cite{matsumura1,matsumura2}) demonstrate that $||u||=||\hat u||\sim e^{-\frac{c}{\nu}t}$ at high frequencies, and this type of decay rate is expected to hold for compressible FENE dumbbell models in periodic domains under zero mean conditions as well.
\end{remark}

Our next result concerns the convergence rate of the solution of (\ref{fene}) under the incompressible limit process $\lambda\rightarrow\infty$. The corresponding limit system is given as follows:

\begin{equation}\label{equ-limit}
    \begin{cases}
        \p_t v + v\cdot\g v - \mu \Delta v + \g\Pi=\div \t \tau,\quad \div v=0,\\
        \p_t \varphi + v\cdot\g \varphi=\divr ( -\g v\cdot R ( \varphi + \pin ) ) + \divr (\pin \gr \frac{\varphi}{\pin}),\\
        \t \tau_{j,k}=\intr (R_j \p_{R_k} \mathcal U)\varphi(x,R,t)\dr,\\
        v(0)=v_0,\quad \varphi(0)=\varphi_0,
    \end{cases}
\end{equation}
where $\Pi$ is the pressure.

We highlight that the convergence rate established in this work exhibits a monotonic decrease as time progresses.
\begin{theorem}\label{thm-limit}
    Let the assumptions in Theorem \ref{thm-exist} hold, with $\varepsilon$ chosen to be smaller if necessary.
    Let $(\elam,\ulam,\plam)$ and $(v,\Pi,\varphi)$ be the solutions of (\ref{fene}) and (\ref{equ-limit}), respectively, and the initial data of (\ref{equ-limit}) satisfies
    $$
    ||v_0||_m \leq \varepsilon,
    $$ \begin{equation}\label{initial-additional}
        \ix \{\P \ulam(0)-v(0)\} \dx =0,\quad
         ||\pru^\nu(0)-v(0)||_{m-1}\leq \varepsilon\nfey, \quad
        \varphi(0) = \plam(0).
    \end{equation}
    Then the following convergence estimate holds:
    \begin{equation*}
        ||\pru^\nu (t) - v (t) ||_{m-1}^2 + ||\psi^\nu (t) -\varphi (t) ||_{m-1,\ler }^2 \leq C b^{-1}(t),
    \end{equation*}
    where
    \begin{equation*}
        b(t)=
        \begin{cases}
            \nu  (1+ \delta t ),\qquad &0<t\leq  \nu,\\
            C_\delta \nu^2\ett,\qquad &t> \nu,
        \end{cases}
    \end{equation*}
\end{theorem}

\begin{remark}
    Combing the convergent results properties in Theorems \ref{thm-exist} and \ref{thm-limit}, we deduce the following convergent estimates:
    $$
    \begin{aligned}
    ||\rho^\nu(t)-1||_m&=||\elam (t)||_m
    \leq C\varepsilon \nfey e^{-\frac{\tc}{\nu}t},\\
        ||\q\ulam (t) ||_m
    &\leq C\varepsilon \nfey  e^{-\frac{\tc}{\nu}t},\\
    ||\P\ulam (t) - v (t) ||_{m-1}
    &\leq ||\pru^\nu (t) - v(t) ||_{m-1} + ||u^\nu(t)||_{m-1} ||\elam(t)||_{m-1}\\
    &\leq C\varepsilon \nfey  e^{-\frac{\tc}{\nu}t}.
    \end{aligned}
    $$
\end{remark}

\begin{remark}
    The initial condition (\ref{initial-additional}) is reasonable and offers a broader range than the usual assumption:
    $$  
    v(0) = \P u(0). 
    $$
    In fact, if the above condition holds, then by the initial data settings in Theorem \ref{thm-exist} and the smallness of $\varepsilon$, we have
    $$
    ||\pru^\nu(0)-v(0)||_{m-1}
    \leq ||\P \ulam(0)-v(0)||_{m-1}
    + ||(\elam\ulam)(0)||_{m-1}\leq \varepsilon\nfey.
    $$
\end{remark}

\subsection{Main difficulties and strategies} \label{sub-difficulties}
In this section, we identify the key challenges addressed in this paper and describe the strategies employed to address them.

When analyzing the global stability of (\ref{fene}), the primary challenge arises from the bad dissipation behavior of the density.
Since we consider the limiting process $\lambda\rightarrow \infty$, and we denote $\nu=2\mu+\lambda$, hence,  this limiting process is equivalently characterized by $\nu\rightarrow\infty$.
First of all, to drive the potential part to zero, we set the energy functional for the potential part as follows:
    $$
    \ep(t)=
     \nu\supt  \ets \bl ||\qulam(s)||_{m}^2+ ||\elam(s)||_{m}^2\brr 
    + \nu^2\int_0^t \ets||\g\qulam(s)||_m^2\d s,
    $$
    where the weight in time characterizes the decay rate.
    
As usual, we investigate the dissipation of the density by the velocity equation in (\ref{fene}):
$$
           \nabla \elam  
           = -\p_t \ulam - \ulam\cdot\g \ulam + \mu\Delta \ulam +( \mu + \lambda )\g\div \ulam +\div \tau^\nu+\ru.
$$
However, the time integrability of the density does not behave well due to the presence of the linear term $( \mu + \lambda )\nabla\div \ulam$ on the right-hand side. 
Indeed, by the $\nu^{-1}$-order of the dissipation of $\nabla\qu$ in $\ep(t)$, we deduce that the dissipation of density exhibits a lower order of $\nu^{-1}$-smallness compared to treating it as energy in $\ep(t)$. That is, while $||\elam||_m^2\sim\nu^{-1}$, we have $\int_0^t ||\nabla \elam||_{m-1}^2 \sim1 $. 
To this end, we set our energy framework of the potential part as:
        $$
    \ee(t)=\int_0^t\ets  ||\g\elam(s)||_{m-1}^2 \d s.
    $$

It is important to emphasize that, compared to the energy of $\elam$ in $\ee(t)$, the dissipation of $\nabla\elam$ in $\ee(t)$ does not benefit from the $\nu^{-1}$-dependent smallness. This presents a significant challenge when estimating the following nonlinear term that arises in the density equation while evaluating $\ee$:
\begin{equation}\label{dif-hard part}
    \begin{aligned}
        &-\nu\ix \gm(\ulam\nabla\cdot\elam)\gm\elam\d x\\
        = &\nu\ix  \div \ulam |\gm\elam|^2 \d x
        -\nu\ix( \gm[\gm,\q \ulam\cdot\nabla]\elam)\gm\elam\d x\\
        &-\nu\ix (\gm[\gm,\P \ulam\cdot\nabla]\elam) \gm\elam\d x.
\end{aligned}
\end{equation}
The first and second terms on the right-hand side can be estimated directly thanks to the strong dissipation of $\g\qu^\nu$. To bound the third term, since the incompressible part $\pu^\nu$ is $\nu^{-1}$-independent, we must treat the entire density $\elam$ as part of the energy in $\ep$ rather than as dissipation in $\ee$ so that the large parameter $\nu$ can be absorbed. 
This implies that we need the following $L^1$-time integrability:
 $$\int_0^t ||\g^s \P \ulam||\d s,\qquad
 1\leq s\leq m.$$
Rather than obtaining $L^1$-time integrability through the energy framework in Besov spaces (see \cite{danchin2017}), we adopt an approach based on estimating the decay properties to establish $L^1$-time integrability.
 To facilitate a more precise decay analysis of $\pu^\nu$, we decompose it into the following two components: 
 $$
 \P u^\nu  = \pru^\nu  -\P(\elam \ulam),
 $$
 Let us briefly explain why we perform the above decomposition.
 The equation of $\P \ulam$ is given by:
 $$
\begin{aligned}
 &\p_t \P\ulam + \P(\ulam\cdot\g \ulam) -\mu\Delta \P\ulam  + \g \P \elam \\
 =&\P\div \tau^\nu  -\nu\frac{\elam}{\rho^\nu}\Delta\qu^\nu
 + \text{other good terms}.
\end{aligned}
 $$
 Among others, due to the existence of the large parameter $\nu$, the nonlinear term $\nu\frac{\elam}{\rho^\nu}\Delta\qu^\nu$ significantly hinders the decay of $\P \ulam$.
 Our strategy is based on the key observation that, in contrast to the equation of $\P\ulam$, the bad nonlinear term $\nu\frac{\elam}{\rho^\nu}\Delta\qu^\nu$ does not appear in the equation of $\pru^\nu$(see (\ref{incompressible})):
 $$
 \p_t\pru +\P \div(\rho u\otimes u)-\mu \P\Delta u = \P \div \tau.
 $$
 As a result, compared to the decay rate of $\P u$, the decay rate of $\pru$ is improved.
 
 For the reason above, it is sufficient to analyze the decay properties of $\pru^\nu$. To this end, we construct the following energy functional:
    $$
    \begin{aligned}
    \ei(t)=&
    \supt a(s) \bl   
    ||\pru^\nu(s)||_{m-1}^2 
    + ||\plam(s)||_{m-1,\ler}^2
    \brr\\
    &+\int_0^t a(s)\bbl  ||\g\pru^\nu(s)||_{m-1
    }^2
    + ||\plam(s)||_{m-1,\hd}^2\bbr \d s,
    \end{aligned}
    $$
    where $a(t)$ is defined in (\ref{at}) and it characterizes the decay rate.
    It is worth emphasizing that achieving this decay rate requires a careful and technical choice of the temporal weight $\ett$, which is motivated by the spectral analysis of the linearized FENE system. 

The rest of our paper is arranged as follows.
In Section \ref{sec-pre}, we introduce the notations and prepare some useful inequalities.
In Section \ref{sec-stability}, we focus on the global stability and the decay estimates of (\ref{fene}), and prove Theorem \ref{thm-exist}.
In Section \ref{sec-incompressible limit}, we utilize the decay property to establish the convergence rates of the incompressible limit, and prove Theorem \ref{thm-limit}.
    
  \section{Preliminaries}\label{sec-pre}
In this section, we introduce the notations and essential lemmas that will be used throughout the paper.
   
    \subsection{Notations}
    In this paper, we will use the following notations.
    We use $f\le g$ to denote $f\leq C g$.
    $\nabla^s$ with a nonnegative integer $s$ denotes the partial derivative $\p_1^{s_1}\p_2^{s_2}\p_3^{s_3}$ with $s_1+s_2+s_3=s$. Also, to simplify the notation, we sometimes abbreviate $f^s:=\nabla^s f$.
     We use $\P$ to denote the Leray projector, defined by $\P f= (\cu I-\nabla\Delta^{-1}\div)f$, and $\q f:= f-\P f$.
    
    Also, we use the abbreviation $B=B(0,1)$. 
    $\p_{R_i}$, $\divr$ and $\gr$ denote the $R_i$ derivative, divergence and gradient in $R$-variable, respectively.
    
    For the norm of Sobolev spaces in $x$-variable, we denote
$$
||f||_m :=||f||_{H^m},\quad
||f|| :=||f||_0,\quad
|||f|||:=||f||_{L^\infty}.
$$

    To characterize the polymer flows, we define the following Hilbert spaces in $R$-variable:
    $$
    \begin{aligned}
    \h L^2 &  := L^2({\dr/\pin })=\left\{
    \psi \ \Big| \ ||\psi||^2_{\h L^2}=\intr |\psi|^2\dri<\infty
    \right\},\\
    \hd &  := \left\{
    \psi \ \Big| \ 
    ||\psi||_\hd^2 = \intr \pin \Big | \gr\frac{\psi}{\pin} \Big |^2 \dr < \infty
\right\}.
    \end{aligned}
    $$

 Next, we define the norms involving $x$ and $R$. For $s\geq 0$,
$$
\begin{aligned}
    ||\psi||_\slerer& :=\sum\limits_{|\alpha|\leq s}\ixr \left|\g^\alpha \psi\right|^2\drx,\\
    ||\psi||_\shyier& :=\sum\limits_{|\alpha|\leq s} \ixr \pin \left | \gr \frac{\g^\alpha \psi}{\pin}\right|^2\dr\dx.
\end{aligned}
$$
Also, we set
$$
||f||_\ler:=||f||_{0,\ler},\quad
||f||_\hyi:=||f||_{0,\hyi}.
$$

We define the linear operator in $R$:
$$
\h L\psi 
 :=\divr ( \gr \psi +\gr\, \mathcal U \psi )
=\divr \left(
\pin \gr \frac{\psi}{\pin}
\right)
$$
with the domain
$$
D(\h L) :=\left\{
\psi\in \h L^2 \ \Big| \
\pin \gr \frac{\psi}{\pin} \in \h L^2, \ \divr \Big (
\pin \gr\frac{\psi}{\pin}\Big )\in \h L^2,\ 
\pin \gr\frac{\psi}{\pin}\cdot n\Big |_{\p B}=0
\right\}.
$$
The boundary condition $\pin \gr\frac{\psi}{\pin}\cdot n\Big |_{\p_B}\!\!\!\!=0$ should be understood in the sense that for any $\phi\in \hd$, 
$$
\intr \phi \h L \psi \dri =-\intr \pin \gr\frac{\phi}{\pin}\gr \frac{\psi}{\pin} \dr .
$$

\subsection{Inequalities}\label{sub-inequality}
The first inequality is about the communicator estimate.
\begin{lemma}\label{le-communicator}
 Suppose $\g f,\g g\in H^{m-1}(\T^3)$, and $m\geq 3$. We have
 $$
 ||[\g^m,f\cdot\g]\,g||\le ||\g f||_{m-1}||\g g||_{m-1}.
 $$
\end{lemma}

Next, we present the following \poin inequalities, which will be utilized frequently throughout this paper.
\begin{lemma}\label{le-poincare}
    Suppose that $(\eta, u, \psi)$ is the solution of (\ref{fene}) satisfying (\ref{mean value 0}), then for any $t\geq 0$, there holds
    $$
    ||\eta (t) || \leq ||\geta(t)||,
    $$
    $$
    ||u(t)||\le ||\g u(t)||,\quad
    ||\qu (t)||\leq ||\g \qu(t)||,
    $$
    $$
    ||\cal M (t)||\leq ||\g \cal M (t)||,\quad
    ||\cal \pru (t)||\leq ||\g \pru (t)||.
    $$
\end{lemma}
For the proof of the above lemma, we refer the readers to \cite{wu2023global}, which builds on the results of \cite{desvillettes2005trend} and Lemma 3.2 in \cite{feireisl2004}.

Finally, we prepare two inequalities to estimate the polymer flow.
The first inequality is the Poincar$\Acute{\text{e}}$ inequality.
  \begin{lemma}[\cite{masmoudi-cpam}]\label{le-poincare-r}
      Assume that $\psi \in \hyi$ and $\intr \psi=0$, then
      \begin{equation*}
          \intr \frac{\psi^2}{\pin} \dr
          \le \int \pin \left| \gr  \frac{\psi}{\pin}\right|^2 \dr.
      \end{equation*}      
  \end{lemma}

 The second inequality is the Hardy type inequality, which is introduced to deal with the singular term $\div \tau$. 
  \begin{lemma}[\cite{masmoudi-cpam}]\label{le-tau}
      Assume that $\psi \in \hyi$ and $\intr \psi=0$, then
\begin{equation*}
          \left(\intr \frac{|\psi|}{1-|R|}\dr \right)^2
      \le \intr \pin 
      \left |
      \gr\frac{\psi}{\pin}
      \right|^2 \dr.
\end{equation*}
  \end{lemma}

\section{Global stability}  \label{sec-stability}
In this Section, we focus on the global stability and decay of system (\ref{fene}). For brevity, we omit the superscript $\nu$. 
Let us begin with the energy framework. We set the basic energy as follows:
    $$
    \begin{aligned}
    \eb(t)=& \supt \bl 
    \ets  || \eta(s)||_{m}^2
    +\ets  ||u (s)||_{m}^2 
    + \ets||\psi(s)||_{m,\h L^2}^2\brr \\
    &    + \int_0^t  \ets \bbl  \mu  ||\g u (s) ||_m^2 + 
    ||\psi(s)||_{m,\hd}^2 \bbr \d s.     
    \end{aligned}
    $$
    Also, we recall the definition of $\ep$, $\ei$ and $\ee$ in Section \ref{sub-difficulties}:
    $$
    \begin{aligned}
    \ep(t)=&
     \nu\supt  \ets \bl ||\qu (s)||_{m}^2+ ||\eta(s)||_{m}^2\brr 
    + \nu^2\int_0^t \ets||\g\qu (s)||_m^2\d s,  \\
    \ei(t)=&
    \supt a(s) \bl   
    ||\pru (s)||_{m-1}^2 
    +   ||\psi(s)||_{m-1,\ler}^2
    \brr\\
    &+\int_0^t a(s)\bbl  ||\g \pru (s)||_{m-1}^2
    + ||\psi (s)||_{m-1,\hd}^2\bbr \d s,\\
    \ee(t)=&\int_0^t\ets  ||\geta (s)||_{m-1}^2 \d s.
    \end{aligned}
    $$
    The main ingredient of the proof of Theorem \ref{thm-exist}  is the establishment of the following energy inequalities.
\begin{prop}\label{prop-key lemma}
        Suppose that the assumptions in Theorem \ref{thm-exist} hold. 
        There exists a constant $\bc$ depends on the constants of  the inequalities in Section \ref{sub-inequality} but independent of $\varepsilon$, $\tc$,  $\nu$, $t$ such that if 
        \begin{equation}\label{eq-assumption}
             \cal E(t) \leq \frac{1}{64\bc^2}, 
        \end{equation} 
        then there exists a positive constant $C$ independent on $\varepsilon$, $\tc$,  $\nu$, $t$ such that
        \begin{align}
                 \eb(t)&\leq C \eb(0)  
    +\frac{ \tc  C}{\nu}  \ee(t) 
    + C \eb^\frac{3}{2}(t)
    + C \ep^\frac{3}{2}(t)
    + C \ee^\frac{3}{2}(t),\label{eb-est}\\   
             \ep(t)&\leq C \ep(0) 
        +C\eb(t) +   \tc  C\ee(t)
        +C\eb^{\frac{3}{2}}(t)
        +C\ep^{\frac{3}{2}}(t)
        +C\ei^{\frac{3}{2}}(t)
        +C\ee^{\frac{3}{2}}(t),\label{ep-est}\\
                 \ei(t)&\leq  C \eb(0) 
         +C\ep(t)
        +C\eb^3(t)
        +C\ep^3(t),\label{ei-est}\\
                 \ee(t)&\leq C\eb(0) 
        +C\eb(t) +  C\ep(t)   
        +C\eb^{\frac{3}{2}}(t)
        +C\ep^{\frac{3}{2}}(t)
        +C\ee^{\frac{3}{2}}(t).\label{ee-est}
        \end{align}
\end{prop}

With the above energy inequalities at hand, we are able to prove Theorem \ref{thm-exist}.  

\begin{proof}[Proof of Theorem \ref{thm-exist}]
    We employ the bootstrap argument (see e.g. P.21 of \cite{tao}). 
    We assume that $\varepsilon\leq \frac{1}{16 \bc}$. Then, there exists a time $t_0>0$ such that 
    $$\cal E(t)\leq 2 \varepsilon^2 
    \leq\frac{1}{64\bc ^2},\quad   \forall 0\leq t\leq t_0.$$ 
    Let $\nu\geq 1$ and $\tc\leq 1$, by Proposition \ref{prop-key lemma}, for $t\leq t_0$, we have
    $$
        \ep(t)  \leq 
        C\cal E(0) 
        + \tc C \ee(t)
        +C\eb^{\frac{3}{2}}(t)
        +C\ep^{\frac{3}{2}}(t)
        +C\ei^{\frac{3}{2}}(t)
        +C\ee^{\frac{3}{2}}(t). 
    $$
    Substituting the above inequality to (\ref{ei-est}) and (\ref{ee-est}), we have
    $$
        \ei(t) + \ee(t) \leq C\cal E(0) 
        + \tc C \ee(t)
        +C\eb^{\frac{3}{2}}(t)
        +C\ep^{\frac{3}{2}}(t)
        +C\ei^{\frac{3}{2}}(t)
        +C\ee^{\frac{3}{2}}(t).
    $$
    Using the above two inequalities, one has
    \begin{equation*}
    \begin{aligned}
        \cal E(t) \leq & C\cal E(0)
        + \tc C \ee(t)
        +C\eb^{\frac{3}{2}}(t)
        +C\ep^{\frac{3}{2}}(t)
        +C\ei^{\frac{3}{2}}(t)
        +C\ee^{\frac{3}{2}}(t)\\
        \leq &C_1\cal E(0)
        + \tc C_1 \cal E(t) + C_1 E^\frac{3}{2}(t),
    \end{aligned}
    \end{equation*}
    where $C_1\geq 1$ is a pure constant.
    Now, we set $\tc=\frac{C_1}{3}$, and make the ansatz that
    \begin{equation*}
        \cal E (t)\leq \frac{1}{36 C_1^2 \bc ^2},
    \end{equation*}
    then, for $\nu\geq \nu_0$, we have
    $$
    \cal E(t)
    \leq C_1\cal E(0) 
    + \frac{1}{3} \cal E(t) 
    + C_1 E^\frac{1}{2}(t) \cal E(t)
    \leq  C_1\cal E(0) + \frac{1}{2} \cal E(t),
    $$
    or 
    $$
     \cal E(t) \leq  2 C_1\cal E(0).
    $$

     If we take the initial data $( \eta_0,u_0,\psi_0 )$ sufficiently small such that
     \begin{equation}\label{ansatz-ini}
     \cal E(0)= \nu ||\eta_0||_m^2 + ||u_0||_m^2 + \nu || \qu_0||_m^2 + ||\psi_0||_{m,\ler}^2
     \leq \frac{1}{144 C_1^3 \bc^2}:=\varepsilon^2,
     \end{equation}
    then 
    \begin{equation}\label{bootstrapbound}
     \cal E (t)\leq \frac{1}{72C_1^2 \bc ^2},
    \end{equation}
    which is exactly half of our ansatz, and the assumption in Proposition \ref{prop-key lemma} automatically holds. 
    Following the bootstrap argument, (\ref{bootstrapbound}) holds for all $t>0$ if the initial data satisfies (\ref{ansatz-ini}). By combining this bound with Lemma \ref{le-poincare}, our proof is accomplished.
\end{proof}

Hence, it remains to prove Proposition \ref{prop-key lemma}. In Section \ref{sec-eb}, we consider the basic energy estimates and prove (\ref{eb-est}).
In  Section \ref{sec-ep}, we introduce the parameter $\nu$ to the potential part of (\ref{fene}) and prove (\ref{ep-est}).
In Section \ref{sec-ei}, we devote to the decay estimates of the incompressible part of the momentum and prove (\ref{ei-est}).
Finally, in Section \ref{sec-ee}, we focus on the dissipation of the density and prove (\ref{ee-est}).

\subsection{Basic estimates}\label{sec-eb}
In this part, we begin our proof with the standard $H^m$-estimates.
\begin{prop}
    Suppose that the assumptions in Proposition \ref{prop-key lemma} hold, then for $t>0$, the following inequality holds:
    $$
    \begin{aligned}
        \eb(t)\leq C \eb(0) 
    +\frac{ \tc  C}{\nu}  \ee(t) 
    + C \eb^\frac{3}{2}(t)
    + C \ep^\frac{3}{2}(t)
    + C \ee^\frac{3}{2}(t).
    \end{aligned}    
    $$
\end{prop}

\begin{proof}
    We begin with the $L^2$-estimate. By integrating (\ref{fene}) with $ (\ett\eta, \ett u ,\ett\psi)$, we have
    \begin{equation}\label{equ-base-l2-est}
    \begin{aligned}
       & \frac{1}{2}\dt \bbl\ett\big(  ||\eta||^2+ ||u||^2+||\psi||_\lerer\big)\bbr\\
        &+\mu \ett||\g u||^2 
        + (\mu + \lambda)\ett ||\div u||^2
        + \ett ||\psi||_\hyier     
        :=\sum\limits_{i=1}^4 I_i,
    \end{aligned}
    \end{equation}
    where
    \begin{equation*}
        \begin{aligned}
            I_1:=&\tc \nu^{-1} \ett\bbl ||\eta||^2+ || u||^2+|| \psi||_\lerer\bbr
        + \ett \ix \div \tau \cdot \qu \dx \\
       &+ \ett \ixr\divr (-\g \qu\cdot R\pin )\psi\drx,\\
     I_2:= & - \ett \ix (u\cdot\g\eta)\eta\dx
        - \ett \ix u\cdot\g u \cdot u\dx \\
        &- \ett \ixr (u\cdot\g \psi)\psi\dri \dx
        - \ett \ix\eta^2\div u\dx\\
        I_3:=& \ett \ix \ru \cdot u\dx ,\\
        I_4:= & \ett \ixr \divr (-\g u\cdot R\psi)\psi\drx.
        \end{aligned}
    \end{equation*}
    Here we use the fact that, by integrating by parts and $-\dfrac{\p_{R_i}\pin}{\pin}=\p_{R_i} \h U$, 
     \begin{equation}\label{linearpart-equal}
     \begin{aligned}
         &\ixr \divr (-\g \pu\cdot R \pin)\psi \drx\\
         =& -\ixr \div \pu\pin\psi  \dr \dx 
         +\sum\limits_{i,j=1}^3 \ixr \p_j \pu_i R_j\p_{R_i}\h U  \psi \dr \dx\\
         = & - \ix \div\tau\cdot \pu \dx.
     \end{aligned}
    \end{equation}

    For linear term $I_1$, we first apply integration by parts along with Lemmas \ref{le-poincare} and \ref{le-tau}, and then choose $\tc$ sufficiently small and $\nu$ sufficiently large to obtain
    \begin{equation*}\label{i1}
    I_1 \leq \ett\bbl 
    \tc \nu^{-1}  ||\geta||^2
    + \frac{1}{2} \ett (\mu+\lambda)||\div u||^2 + \frac{1}{2}  \mu || \g u||^2+ \frac{1}{2} ||\psi||_\hyier \bbr.
    \end{equation*}
    
    For $I_2$, using integrating by parts, $\div \pu=0$, Lemmas \ref{le-poincare}, \ref{le-poincare-r} and Sobolev embedding, we have
    \begin{equation*}\label{i3}
    \begin{aligned}
    I_2
    =  & \ett \ix (\qu\cdot\geta)\eta\dx
    - \ett \ix (\qu\cdot\g u)\cdot u \dx\\
    & - \ett \ixr (\qu \cdot\g \psi) \psi\drx \\
    \le &  \ett ||\g\qu||\,||\geta||\,|||\eta|||
    +  \ett ||\g\qu||\,||\g u||\,|||u|||\\
    & + \ett |||\qu||| ||\psi||_{1,\ler}^2\\
    \le & \bl ||u||_2+||\eta||_2\brr  \ett \bbl
      ||\g u||^2 + ||\geta||^2 +  
 ||\psi||_{1,\hyi}^2 \bbr.
    \end{aligned}
   \end{equation*}

    For $I_3$ and $I_4$, by using integrating by parts and Lemmas \ref{le-poincare}, \ref{le-tau} and assumption (\ref{eq-assumption}), one has
    \begin{equation*}\label{i7}
    \begin{aligned}
    I_3=&\mu \ett  \ix  (\g\er\cdot
    \g u) \cdot u \dx
    +\mu \ett  \ix \er|\g u|^2 \dx\\
    &+( \mu + \lambda ) \ett \ix \div u (\g\er\cdot u)\dx \\
    &+( \mu + \lambda ) \ett \ix \er |\div u|^2 \dx
    -  \ett \ix\eta(\geta\cdot\qu)\dx\\ 
    &+ \ett \ix(\g\er \cdot\tau)\cdot u \dx
    + \ett \ix \er(\tau:\g u) \dx\\
    \le &
    \mu  \ett |||u|||\,||\geta||\,||\g u|| 
    + \mu  \ett |||\eta|||\,||\g u||^2\\
    &+   \ett |||u|||\,||\geta||\nu||\g\qu||
    + \ett  |||\eta|||\,\nu||\g\qu||^2\\
    &+  \ett |||\eta|||\,||\geta||\,||\g\qu ||
    + \ett |||u|||\,||\geta||\,||\psi||_\hyier\\
    &+ \ett |||\eta|||\,||\psi||_\hyi||\g u|| \\
    \le & \bl  ||u||_2 + ||\eta||_2 \brr 
    \ett\bbl  ||\geta||^2 + ||\g u||^2 
    + \ne||\g\qu||^2
    +||\psi||_\hyier \bbr
    \end{aligned}
    \end{equation*}
and 
$$
\begin{aligned}
        I_4=& \ett \ixr \Big( \g u\cdot R\pin^{-\frac{1}{2}}\psi \Big)
    \Big(\pin^{\frac{1}{2}}\gr \frac{\psi}{\pin}  \Big)  \dr\dx\\
    \le & ||u||_1\,  \ett ||\psi||_{1,\hyi}^2.
\end{aligned}
$$
    
    Substituting the above estimates of $I_1$-$I_4$ into (\ref{equ-base-l2-est}), and integrating in time, we obtain that
    \begin{equation}\label{est-base-l2}
        \begin{aligned}
            & \ett ||\eta(t)||^2+ \ett ||u(t)||^2
            + \ett ||\psi(t)||_\lerer\\
           & +\mu  \it\ets ||\g u||^2 \d s 
             +(\mu+\lambda) \it \ets ||\div u ||^2 \d s 
            + \it\ets ||\psi||_\hyier \d s \\
            \leq & ||\eta(0)||^2+||u(0)||^2+||\psi(0)||_\lerer
            +\tc \nu^{-1}\it  \ets ||\geta||^2\d s \\
            & +C \supt \bl ||u||_2 + ||\eta||_2 \brr
    \it \ets \Big\{ ||\geta||^2 + ||\g u||^2 + \nu^2 ||\g\qu ||^2 + ||\psi||_{1,\hyi}^2\Big\}\d s.\\
        \end{aligned}
    \end{equation}

    Next, we turn to the $\dot H^m$-estimate.
    By applying $\gm$ to (\ref{fene}) and multiplying the resulting equation with $\ett(\gm\eta,\gm u,\gm\psi)$, we obtain that
    \begin{equation}\label{equ-base-hm-est}
    \begin{aligned}
       & \frac{1}{2}\dt \Big\{\ett\big( ||\gm\eta||^2+ ||\gm u||^2+||\gm \psi||_\lerer\big)\Big\}\\
       & +\mu\ett ||\gm\g u||^2 
        + (\mu + \lambda)\ett ||\gm \div u||^2
        +\ett ||\gm \psi||_\hyier
        =\sum\limits_{i=1}^4 J_i,
    \end{aligned}
    \end{equation}    
    where
    $$
    \begin{aligned}
        J_1:= & \tc \nu^{-1} \ett\bbl ||\gm\eta||^2+ ||\gm u||^2+||\gm \psi||^2\bbr\\
        &+ \ett \ix \gm\div \tau \cdot \gm\qu \dx\\
        &+\ett \ixr \gm \divr (-\g\qu\cdot R\pin )\gm\psi\drx,\\
        J_2:= & -\ett \ix \gm(u\cdot\g\eta)\gm\eta\dx
        -\ett \ix \gm(u\cdot\g u) \cdot \gm u\dx \\
       & -\ett \ixr \gm(u\cdot\g \psi)\gm\psi\dri \dx,\\
       J_3:= & -\ett \ix\gm(\eta\div u)\gm\eta\dx\\
        &+\ett \ixr \gm\divr (-\g u\cdot R\psi)\gm\psi\drx,\\
        J_4:= & \ett \ix \gm\ru \cdot \gm u\dx.
    \end{aligned}
    $$
    
    For linear term $J_1$, by choosing proper $\tc $ and $\nu$, we directly have
    $$
     J_1
    \leq 
    \ett\Big\{ \tc \nu^{-1}  ||\gm\eta||^2
    +\frac{\mu}{2} ||\gm \g u||^2
    +\frac{1}{2} (\mu+\lambda)||\gm\div u||^2
    +\frac{1}{2}||\gm \psi||_\hyier\Big\}.
    $$
    
    For the convection term $J_2$, we use integrating by parts and Lemma \ref{le-communicator} to obtain 
    $$
    \begin{aligned}
        J_2
       = &\ett \ix \div u \Big\{|\gm\eta|^2+ |\gm u|^2\Big\} \dx 
        + \ett \ixr \div u |\gm\psi|^2\drx\\
        &+ \ett \ix ([\gm, u\cdot\g]\eta) \gm\eta \dx
        + \ett \ix ([\gm, u\cdot\g]u)\cdot \gm u \dx\\
        & + \ett \ix ([\gm, u\cdot\g]\psi)\gm\psi \drx\\
        \le & ||\g u||_{m-1} \ett \bbl ||\g\eta||_{m-1}^2 +||\g u||_{m-1}^2 + ||\g\psi||_\mmhyier \bbr.
    \end{aligned}
    $$
    
    For $J_3$, by integrating by parts,
    $$
    \begin{aligned}
    J_3 &\le 
    \ett ||\eta||_{m}||\g u||_m||\gm\eta||
    + \ett ||\psi||_\mler||\g u||_{m}||\gm\psi||_\hyi\\
    &\le \bl ||\eta||_m+||\psi||_\mler\brr 
    \ett \bbl  ||\g u||_m^2 +  ||\gm \eta||^2 + ||\gm\psi||_\hyier \bbr.
    \end{aligned}
    $$
    
      By the assumption (\ref{eq-assumption}), Lemma \ref{le-poincare} and $m\geq 3$, there exists a polynomial $\cal P$ such that for $1\leq s\leq m$, there holds
      \begin{equation}\label{low bound}
      ||\nabla\big(\er\big)||_{s-1}\leq 
      ||\er||_s
      \leq 
      \cal P(|||\eta|||_{[\frac{s}{2}]-1}, ||\g^{[\frac{s}{2}] }\eta||_{L^4})||\geta||_{s-1}\le ||\nabla\eta||_{m-1}.
      \end{equation}
      Using the above property and integrating by parts, we can bound $J_4$ by
      $$
      \begin{aligned}
          J_4\le &  \ett|| \gmm \bl\eta\g \eta ,\,\mu  \er\Delta u ,\,\nu  \er\g\div u ,\, \er\div \tau \brr || \, ||\g\gm u||\\
          \le & ||\eta||_{m}
          \ett \bbl  ||\geta||_{m-1}^2 + ||\g u||_m^2 + \ne||\g \qu||_m^2 +||\g \psi||_\mmhyier \bbr .      \end{aligned}
      $$

      Gathering the estimates of $J_1$-$J_4$ together, substituting them into (\ref{equ-base-hm-est}) and integrating the resulting inequality in time, we deduce
    \begin{equation}\label{est-base-hm}
        \begin{aligned}            
        &\ett||\gm\eta(t)||^2
        +\ett||\gm u(t)||^2
        +\ett||\gm \psi(t)||_\lerer\\
            &+\mu \it\ets||\g\gm u||^2 \d s 
            +(\mu+\lambda) \it\ets||\gm\div u||^2 \d s 
            + \it\ets||\gm \psi||_\hyier \d s \\
            \leq &
            ||\gm\eta(0)||^2
        +||\gm u(0)||^2
        +||\gm \psi(0)||_\lerer
        + \frac{ \tc C }{\nu} \it \ets ||\geta ||_{m-1}^2 \d s\\
            &+ C \supt \bl ||u(s)||_m + ||\eta(s)||_m + ||\psi(s)||_\mler\brr\\
   &\qquad\times  \it \ets\Big\{  ||\geta||_{m-1}^2 
    + ||\g u||_{m}^2 
    + \nu^2 ||\g\qu ||_{m}^2 
    + ||\g\psi||_\mmhyier\Big\}\d s.
        \end{aligned}
    \end{equation}      
    
      Combining (\ref{est-base-hm}) with (\ref{est-base-l2}), and using Lemma \ref{le-poincare} and the norm equivalence 
      \begin{equation}\label{equ-normequi}
          ||f||_{  H^{m} }\sim || f|| + ||\gm f||,
      \end{equation}
      yields
      $$
        \begin{aligned}            
        &\supt\bl \ets||\eta(s)||_{m}^2
        +\ets|| u(s)||_{m}^2
        +\ets||\psi(s)||_\mlerer\brr \\
            &+\mu \it\ets||\g u||_m^2 \d s 
            +\nu  \it\ets||\div u||_m^2 \d s 
            + \ets|| \psi||_\mhyier \d s \\
            \leq & 
            C ||\eta(0)||_m^2 
            + C ||u(0)||_m^2 
            +C ||\eta(0)||_\mlerer
            +\frac{ \tc  C}{\nu} \it \ets ||\geta ||_{m-1}^2 \d s\\
            &+ C \supt \bl ||u(s)||_m + ||\eta(s)||_m +  ||\psi(s)||_\mler\brr\\
   &\qquad\times  \it \ets\Big\{  ||\geta||_{m-1}^2 
    + ||\g u||_{m}^2 
    + \nu^2 ||\g\qu ||_{m}^2 
    + ||\g\psi||_\mmhyier\Big\}\d s\\
    \leq & C \eb(0)  
    +\frac{ \tc  C}{\nu}\ee(t)
    + C \eb^\frac{3}{2}(t)
    + C \ep^\frac{3}{2}(t)
    + C \ee^\frac{3}{2}(t),
        \end{aligned}      
      $$
    which is exactly (\ref{eb-est}).
\end{proof}

\subsection{Estimates for potential part}\label{sec-ep}
In this section, we focus on estimating the potential part of (\ref{fene}) and aim to establish (\ref{ep-est}) and (\ref{ee-est}).
As outlined in the introduction, the primary challenge lies in the lack of $\nu^{-1}$-dependent smallness of the dissipation of density. To overcome this problem, we take advantage of the damping effect in $\ei$ to obtain the enhanced time-integrability of the incompressible part of velocity.

We first deduce the system of the potential part of (\ref{fene}). By applying projector $\q$ to the second equation of (\ref{fene}), we obtain 
\begin{equation}\label{eq-compressible}
    \begin{cases}
        \p_t \eta + u\cdot\geta + \div\qu = -\eta\div \qu,\\
        \p_t \qu + \q ( u\cdot\nabla \qu ) +\geta -\nu\Delta\qu =\rqu,
    \end{cases}
\end{equation}
where
\begin{equation*}\label{rqu-def}
    \rqu = \q(\frac{1}{\rho}\div\tau )-\q ( u\cdot\nabla\pu  )- \mu\q( \er\Delta u )
-(\mu+\lambda)\q ( \er\nabla\div u ) -\q (\eta\geta).
\end{equation*}

The following result proves (\ref{ep-est}).
\begin{prop} \label{prop-e2}
    Suppose that the assumptions in Proposition \ref{prop-key lemma} hold, then for $t>0$, the following inequality holds:
    $$
    \ep(t)\leq C\ep(0) + C\eb(t) + \tc C \ee(t) 
        + C\eb^\frac{3}{2}(t)
        + C\ep^\frac{3}{2}(t)
        + C\ei^\frac{3}{2}(t)
        + C\ee^\frac{3}{2}(t).
    $$
\end{prop}
\begin{proof}
Since both $\eta$ and $\qu$ satisfies the \poin inequality, we only need to consider the $\dot H^m$-estimate. Applying $\gm$ to (\ref{eq-compressible}) and taking the inner product of $(\nu\ett\gm\eta,\nu\ett\gm\qu)$ with the resulting equations, we have
\begin{equation}\label{equ-ep-hm}
     \begin{aligned}
         \frac{1}{2}\dt \bbl \nu\ett \big ( ||\gm\eta||^2
         +||\gm\qu||^2 \big ) \bbr
         +\ne\ett||\g\gm\qu||^2  
         :=\sum\limits_{i=1}^{5} K_i.
     \end{aligned}
\end{equation}
where 
$$
\begin{aligned}
    K_1:= & \tc \ett \bbl ||\gm\eta||^2+||\gm\qu||^2 \bbr,\\
    K_2:= & -\nu\ett\ix \gm( u\cdot\g\eta )\gm\eta\dx,\\
    K_3:= & -\nu\ett\ix \gm( u\cdot\g\qu )\cdot\gm\qu\dx,\\
    K_4:= & -\nu\ett\ix\gm( \eta\div u )\gm\eta \dx,\\
    K_5:= & \nu\ett\ix\gm\rqu\cdot\gm\qu\dx.
\end{aligned}
$$

For $K_1$, we use \poin inequality  and set $\tc\leq \frac{\nu^2}{4}$ to obtain
\begin{equation}\label{k1}
    K_1\leq \tc  \ett ||\gm\eta||^2 +\frac{1}{4}\nu^2 \ett  ||\g \gm\qu||^2 .
\end{equation}

For $K_2$, we separate it into two terms:
$$
\begin{aligned}
K_2=&-\nu\ett\ix (u\cdot\g\gm\eta)\gm\eta\dx\\
&- \nu\ett\ix ([\gm,u\cdot\g]\eta)\gm\eta\dx
:=K_{21}+K_{22}.
\end{aligned}
$$
The bound of $K_{21}$ is straight forward:
$$
\begin{aligned}
K_{21}=&\nu\ett\ix \div u|\gm\eta|^2\dx\\
\le&||\gm\eta||\ett 
\bbl \ne||\g\qu||_2^2
+   ||\gm\eta||^2 \bbr
\end{aligned}
$$
To control $K_{22}$, we need to take advantage of the decay property of $\pru$. 
Since 
\begin{equation}\label{u-decompose}
    u=\pru - \P \eu + \qu,
\end{equation}
we can use Lemma \ref{le-communicator} and Lemma \ref{le-poincare} to obtain
   $$
   \begin{aligned}
        K_{22}=&-\nu\ett\ix ([\gm,\pru\cdot\g]\eta)\gm\eta\dx \\
        &+\nu\ett\ix ([\gm,\P(\eta u )\cdot\g]\eta)\gm\eta\dx\\
        &-\nu\ett\ix ([\gm,\qu\cdot\g]\eta)\gm\eta\dx\\
        \le & \nu \ett ||\geta||_{m-1}^2 \bl ||\g\pru||_{m-1}
        + ||\g\eu||_{m-1}
        +||\g\qu||_{m-1}\brr\\
        \le & \nu  ||\geta||_{m-1}^2 
        \ett \bbl ||\g\pru||_{m-1} + ||\geta||_{m-1}^2 + ||\g u||_{m-1}^2\bbr \\
        & + ||\eta||_m 
        \ett \bbl ||\geta||_{m-1}^2 + \nu^2\ett||\g\qu||_{m-1}^2\bbr. 
   \end{aligned}
$$
Hence, $K_2$ can be bounded by
\begin{equation}\label{k2}
    \begin{aligned}
        K_2\le &\bl
      ||\eta||_m
        +  \nu||\geta||_{m-1}^2\brr
        \ett \bbl 
         ||\g u||_{m-1}^2
        +\ne||\g\qu||_{m-1}^2 
        +||\geta||_{m-1}^2\bbr\\
        &+ \nu\ett||\geta||_{m-1}^2||\g\pru||_{m-1}.
    \end{aligned}
\end{equation}

For $K_3$, we use integrating by parts and Lemma \ref{le-communicator}
to obtain
\begin{equation}\label{k3}
    \begin{aligned}
        K_3=&\nu\ett\ix \div u |\gm\qu|^2\dx\\
        & -\nu\ett\ix ([\gm,u\cdot\g]\qu)\cdot\gm\qu\dx\\
        \le & \nu\ett||\g u||_2||\gm\qu||^2
        + \nu\ett||\g u||_{m-1}||\g\qu||_{m-1}||\gm\qu||\\
        \le & ||\g u||_{m-1}\nu\ett||\g\qu||_{m-1}^2.
    \end{aligned}
\end{equation}

Direct computation shows
\begin{equation}\label{k4}
\begin{aligned}
    K_4 & \le \ett ||\eta||_m||\g\qu||_m||\gm\eta||\\
   & \le  ||\eta||_m
    \ett \bbl ||\geta||_{m-1}^2 + \ne ||\g\qu||_m^2\bbr .
\end{aligned}
\end{equation}

 Using (\ref{eq-assumption}) and (\ref{low bound}), we can bound $K_5$ by
\begin{equation}\label{k5}
\begin{aligned}
    K_{5}\leq & \nu\ett
    ||
    \mu\er\Delta u,\,
    \nu \er\g\div u,\,
    \eta\geta,\,
    u\cdot\g\pu,\,
    \frac{1}{\rho}\div\tau||_{\dot H^{m-1}}\,
    ||\gmp\qu||\\
    \leq &  C\nu\ett ||\eta||_m
     \bl ||\g u||_m + \nu ||\g\qu||_m 
     +  ||\geta||_{m-1} \brr ||\gmp\qu||\\
     & + C\nu\ett||u||_{m-1}||\g u||_{m}||\gmp\qu||\\
    &  + C\nu\ett ||\g\psi||_\mmhyi||\gmp\qu||\\
     \leq  & C \bl ||u||_m + ||\eta||_m \brr 
     \ett \bbl  ||\g u||_m^2 
     + \nu^2 ||\g\qu||_m^2 
     + ||\geta||_{m-1}^2\bbr\\
    & + C\ett||\g\psi||_{{m-1},\hyi}^2 + \frac{1}{4}\ne\ett||\gmp\qu||^2.
\end{aligned}
\end{equation}

Now, plugging (\ref{k1})-(\ref{k5}) into (\ref{equ-ep-hm}) yields
\begin{equation*}
    \begin{aligned}
        & \dt \bbl \nu\ett  \big(||\gm\eta||^2 
        +  ||\gm\qu||^2 \big)\bbr + \ne\ett||\gmp\qu||^2\\
        \leq & \tc C \ett||\gm\eta||^2 
        + C\ett||\g\psi||_\mmhyier
         +\nu\ett||\geta||_{m-1}^2  ||\g\pru||_{m-1}\\
        &+\bbl ||u||_m + ||\eta||_{m} 
        + \nu||\eta||_{m}^2 \bbr
        \ett\bbl 
         ||\g u||_m^2
        +||\geta||_{m-1}^2 
        + \ne||\g\qu||_m^2
         \bbr.
    \end{aligned}  
\end{equation*}

      Integrating the above inequality in time, and applying Poincar$\Acute{\text{e}}$ inequality and the norm equivalence (\ref{equ-normequi})
      to the resulting inequality, we have
\begin{equation}\label{est-qu-final}
    \begin{aligned}
        & 
         \supt \bbl \nu\ets ||\eta(s)||_{m}^2
        +  \nu\ets ||\qu(s)||_{m}^2\bbr 
         + \ne \ix\ets ||\g\qu(s)||_m^2 \d s\\
        \leq &  C  \nu||\eta(0)||_m^2 +  C \nu||\qu(0)||_m^2 \\
        & + C\it \ets\bbl||\g\psi(s)||_\mmhyier 
        +\tc ||\geta(s)||_{m-1}^2 \bbr \d s\\
        & + C\supt\bbl
          ||u(s)||_m 
         + ||\eta(s)||_{m-1}
         + \nu||\geta(s)||_{m-1}^2
        \bbr\\
        &\qquad\times
        \it \ets
        \bbl  ||\g u(s)||_m^2
        +\ne||\g\qu(s)||_m^2
        +||\geta(s)||_{m-1}^2\bbr\d s\\
        &+ C \supt \bbl \nu\ets||\eta(s)||_{m}^2\bbr
        \it ||\g\pru(s)||_{m-1}\d s\\
        \leq & C\ep(0) + C\eb(t) + \tc C \ee(t) 
        + C\eb^\frac{3}{2}(t)
        + C\ep^\frac{3}{2}(t)
        + C\ei^\frac{3}{2}(t)
        + C\ee^\frac{3}{2}(t),
    \end{aligned}
\end{equation}
where in the last inequality, we use assumption (\ref{eq-assumption}) and
$$
\begin{aligned}
&\it ||\g\pru(s)||_{m-1}\d s \\
\le & \bbl\it a(s)||\g\pru(s)||_{m-1}^2\d s\bbr^\frac{1}{2}
\bbl\it a^{-1}(s)\d s\bbr^\frac{1}{2}
\le \ei^\frac{1}{2}(t).
\end{aligned}
$$
Combining (\ref{est-qu-final}) with the definition of $\ep$, our proof is acomplished.
 \end{proof}
 \subsection{Estimates for incompressible part}\label{sec-ei}
 In this section, we consider the decay estimates of the incompressible part of momentum and prove (\ref{ei-est}).
 For convenience, we rewrite the incompressible part of (\ref{fene}) as follows
 \begin{equation}\label{incompressible}
     \begin{cases}
         \p_t\pru +\P \div(\rho u\otimes u)-\mu \P\Delta u = \P \div \tau,\\
         \psi_t+ u\cdot\g \psi = \divr (-\g u\cdot R( \psi+\pin )) + \h L\psi.
     \end{cases}
 \end{equation}

 \begin{prop}
    Suppose that the assumptions in Proposition \ref{prop-key lemma} hold, then for $t>0$, the following inequality holds:
    $$
    \begin{aligned}
        \ei(t)\leq  C \eb(0) 
         +C\ep(t)
        +C\eb^3(t)
        +C\ep^3(t).
    \end{aligned}
    $$     
 \end{prop}
\begin{proof}
    We begin with the $L^2$-estimates. By dotting (\ref{incompressible}) with $(\pru,\psi)$ and using (\ref{linearpart-equal}), we find
    \begin{equation}\label{equ-incompressible-l2}
    \begin{aligned}
        \frac{1}{2}\frac{\d}{\d t}\bbl  ||\pru||^2 + ||\psi||_\lerer\bbr
        + \mu||\g\pru||^2 + ||\psi||_\hyier
        = \sum_{i=1}^4 M_i.
    \end{aligned}
    \end{equation}
    where
    $$
    \begin{aligned}
        M_1:= & \mu\ix \g\P \eu \cdot\g\pru \dx,\\ 
        M_2:= & \ixr \divr ( \g \eu\cdot R \pin)\psi\drx,\\
        M_3:= & -\ix \div (\rho u\otimes u)\cdot\pru\dx,\\
        M_4:= & -\ixr (u\cdot\g\psi) \psi\drx\\
        &+\ixr \divr (-\g u\cdot R\psi)\psi\drx\\
       & +\ixr \divr (-\g \qu\cdot R\pin)\psi\drx.
    \end{aligned}
    $$

    To bound $M_1$, $M_2$ and $M_3$, we first estimate $||\g\eu||$. From now on, we use $\bc\geq 1$ to represent a constant that varies depending on the case and is related to $\mu$, Sobolev embedding, \poin inequality. 
    By (\ref{u-decompose}), we have
    \begin{equation}\label{eu}
    ||\eta u||_1\leq \bc ||\eta||_2||\g\pru||
    + C ||\eta||_2^2||\g u ||
    + C ||\eta||_2||\g\qu||.
    \end{equation}
    Hence, using (\ref{eq-assumption}) and Lemma \ref{le-poincare},
    we can bound $M_1$, $M_2$ and $M_3$ by
    \begin{equation*}\label{m1}
    \begin{aligned}
        M_1\leq&  \bc \mu ||\eta||_2||\g\pru||^2  + \frac{1}{8}\mu ||\g\pru||^2
        + C ||\eta||_2^4 ||\g u||^2
        +C ||\eta||_2^3 ||\g \qu||^2\\
        \leq & \frac{3}{16} \mu   ||\g\pru||^2
        + C ||\eta||_2^4 ||\g u||^2
        +C ||\eta||_2^3 ||\g \qu||^2,\\
            M_2\leq &\bc||\eta||_2||\g\pru||\,||\psi||_\hyi
            + C||\eta||_2^2||\g u||\,||\psi||_\hyi
            + C||\eta||_2||\g \qu||\,||\psi||_\hyi\\
            \leq & \frac{1}{16}\mu ||\g\pru||^2 + \frac{1}{8}||\psi||_\hyier
            + C ||\eta||_2^4 ||\g u||^2
        +C ||\eta||_2^2 ||\g \qu||^2,\\
            M_3=&
            \ix (\pru\otimes u)\cdot\g\pru\dx 
            - \ix (\P\eu\otimes u)\cdot\g\pru\dx\\
            &+\ix (\qu\otimes u)\cdot\g\pru\dx\\
            \leq & \bc ||u||_1||\g\pru||^2
            + C ||\eta u||_1||u||_1||\g\pru||
            + C ||u\|_1||\g\qu||\,||\g\pru||\\
            \leq &\frac{3}{16}\mu ||\g\pru||^2 
           + C ||\eta||_2^4 ||\g u||^2
        +C ||\eta||_2^2 ||\g \qu||^2. 
    \end{aligned}
    \end{equation*}
    
    By choosing $\nu\geq 16\bc $, we can use (\ref{eq-assumption}) to bound $M_4$ as follows:
    \begin{equation*}\label{m4}
        M_4\leq \frac{1}{4}||\psi||_\hyier + C||\g\qu||^2.
    \end{equation*}

    Substituting the estimates of $M_1$-$M_4$ into (\ref{equ-incompressible-l2}) yields
    \begin{equation}\label{est-incompressible-l2}
    \begin{aligned}
        &\frac{\d}{\d t}\bbl  ||\pru||^2 + ||\psi||_\lerer\bbr
        + \mu||\g\pru||^2 + \frac{5}{4} ||\psi||_\hyier\\
        \leq & C ||\eta||_2^4||\g u||^2 
        +C ||\g\qu||^2.
    \end{aligned}
    \end{equation}    

    Next, we consider the $\dot H^{m-1}$-estimates.
    By applying $\gmm$ to (\ref{incompressible}), integrating the resulting equation with ($\gmm\pru,\gmm\psi$) and utilizing (\ref{linearpart-equal}), we have
    \begin{equation}\label{equ-incompressible-hmm}
        \begin{aligned}
        \frac{1}{2}\frac{\d }{\d t} \bbl  ||\pru\mm||^2 + ||\psi\mm||_\lerer \bbr
        +\mu ||\g\pru\mm||^2 +||\psi\mm||_\hyier
        =  \sum_{i=1}^6 N_i,
        \end{aligned}
    \end{equation}
    where
    $$
    \begin{aligned}
        N_1 := &  \mu\ix \g \gmm \P \eu \cdot\g  \pru\mm \dx,\\
        N_2 := & \ixr \divr ( \g  \gmm  \eu\cdot R \pin) \psi\mm\drx,\\
        N_3 := & -\ix \gmm \div (\rho u\otimes u)\cdot\pru\mm\dx,\\
        N_4:= & -\ixr \gmm(u\cdot\g\psi) \psi\mm\drx,\\
        N_5 := & \ixr \gmm\divr (-\g u\cdot R\psi)\psi\mm\drx,\\
        N_6 := & \ixr \divr (-\g \qu\mm\cdot R\pin)\psi\mm\drx.
    \end{aligned}
    $$

        To bound $N_1$, $N_2$ and $N_3$, we first estimate $||\gm\eu||$. By Lemma \ref{le-poincare}, for $s\leq m$, there holds
    \begin{equation}\label{eu-hs}
    \begin{aligned}
    ||\g^s\eu||
    \leq &  C ||\gm\eta||( ||\gm\pru||+ ||\gm\eu||+ ||\gm\qu|| )\\
    \leq &\bc||\eta||_m||\g\pru\mm|| + C ||\eta||_m^2 ||\g u||_{m-1} 
    + C ||\eta||_m||\g\qu||_{m-1}.
    \end{aligned}
    \end{equation}
    Combining the above inequality with Lemma \ref{le-poincare} and assumption (\ref{eq-assumption}), we can bound $N_1$, $N_2$ and $N_3$ by
    \begin{equation*}\label{n1}
        \begin{aligned}
            N_1
            \leq & \bc \mu||\eta||_m\,||\g\pru\mm||^2
            + C||\eta||_m^2||\g u||_{m-1}\,||\g\pru\mm||\\
            & + C||\eta||_m||\g\qu||_{m-1}\,||\g\pru\mm||\\
            \leq & \frac{3}{16}\mu ||\g\pru\mm||^2
             + C||\eta||_m^4||\g u||_{m-1}^2 
             + C ||\eta||_m^2  ||\g\qu ||_{m-1}^2,\\
            N_2=& \ixr \gmm( \g\P\eu\cdot R\pin )\gr\frac{\psi\mm}{\pin}\drx\\
            \leq & \frac{1}{16}\mu||\g\pru\mm||^2
            + \frac{1}{8}||\psi\mm||_\hyier\\
            & + C||\eta||_m^4||\g u||_{m-1}^2 
              + C ||\eta||_m^2  ||\g\qu ||_{m-1}^2,\\
            N_3\leq & \bc ( ||\pru||_{m-1}
            +||\q\eu||_{m-1}
            +||\qu||_{m-1})||u||_{m-1}||\g\pru\mm||\\
            \leq  & \frac{3}{16}\mu ||\g\pru\mm||^2 
            + C||\eta||_m^4||\g u||_{m-1}^2 
             + C ||\eta||_m^2 ||\g\qu ||_{m-1}^2.
        \end{aligned}
    \end{equation*}
    
    For $N_4$, $N_5$ and $N_6$, using Lemma \ref{le-poincare}, we have
    \begin{equation*}\label{n4n6}
        \begin{aligned}
            N_4 = & \ixr \div u|\psi\mm|^2\drx
            -\ixr ([\gm,u\cdot\g]\psi)\psi\mm\drx\\
            \leq &  \bc ||u||_m||\psi\mm||_\hyier
            \leq  \frac{1}{8}||\psi\mm||_\hyier,\\
        N_5\leq & \bc ||\g u||_{m-1}( ||\psi||_\hyi + ||\psi\mm||_\hyi )||\psi\mm||_\hyi\\
        \leq &\frac{1}{16}||\psi||_\hyier
        + \frac{3}{16}||\psi\mm||_\hyier, \\
        N_6\leq & \frac{1}{16}||\psi\mm||_\hyier + C||\g\qu||_{m-1}^2.
        \end{aligned}
    \end{equation*}

    Gathering the bounds of $N_1$-$N_6$ together and substituting them into (\ref{equ-incompressible-hmm}) yields
    \begin{equation*}
        \begin{aligned}
            &\frac{\d }{\d t} \bbl  ||\pru\mm||^2 + ||\psi\mm||_\lerer \bbr
        +\mu ||\g\pru\mm||^2 +||\psi\mm||_\hyier\\ 
        \leq &\frac{1}{8}||\psi||_\hyier
        + C ||\eta||_m^4||\g u||_{m-1}^2
        + C ||\g\qu||_{m-1}^2.
        \end{aligned}
    \end{equation*}
    By combining the above inequality with (\ref{est-incompressible-l2}), we deduce
    \begin{equation*}
        \begin{aligned}
            & h'(t) 
            +g(t)
            \leq  C ||\eta||_m^4||\g u||_{m-1}^2
        + C ||\g\qu||_{m-1}^2,
        \end{aligned}
    \end{equation*}
    where
    $$
    \begin{aligned}
    h(t):=&||\pru(t)||^2+||\pru\mm(t)||^2 +||\psi(t)||_\lerer+ ||\psi\mm(t)||_\lerer,\\
    g(t):= &\mu ||\g\pru(t)||^2 
            + \mu ||\g\pru\mm(t)||^2
            + ||\psi(t)||_\hyier
            +  ||\psi\mm(t)||_\hyier.
    \end{aligned}
    $$
    Also, we recall the definition of the continuous function $a(t)$:
    \begin{equation*}
        a(t)=
        \begin{cases}
            (1+\delta t )^2, \qquad &0< t \leq \nu,\\
            C_\delta \nu^2\ett,\qquad &t> \nu,
        \end{cases}
    \end{equation*}
    We set $C_\delta:=\frac{(1+\delta \nu )^2}{\nu^2e^{2\tc}}\sim 1 $ such that $a(t)$ is continuous, and for we have $t\leq \nu$, $a(t)\le \nu^2\ett$.
    Also, we assume $0<\tc\leq 1$.
    By Lemma \ref{le-poincare}, \ref{le-poincare-r}, there exists a constant $C_P\geq 1$ such that $ h(t)\leq C_P g(t)$. Hence, by choosing 
    \begin{equation}\label{delta}
    \nu\geq 4 C_P 
    \quad\text{and} \quad 
    \delta\leq \frac{\mu }{4 C_P}, 
    \end{equation}
    we have
    $$
    a'(t)h(t)\leq \frac{1}{2}a(t) g(t),
    $$
    which implies
    $$
    \begin{aligned}
    &\dt\bl a(t) h(t)\brr + \frac{1}{2}a(t) g(t)\\
    \leq & C a(t)||\eta||_m^4||\g u||_{m-1}^2
        + C  a(t)||\g\qu||_{m-1}^2.
    \end{aligned}
    $$
    By integrating the above inequality in time, we find
    $$
    \begin{aligned}
        &a(t)h(t) + \it a(s) g(s)\d s \\
        \leq & h(0) 
        + C \nu^2 \it \ets||\g\qu||_{m-1}^2\d s\\
        & +  C \bbl  \nu^2 \supt \ets||\eta(s)||_m^4 \bbr
        \it ||\g u||_{m-1}^2\d s\\
        \leq & \eb(0) + C \ep(t) + C \eb^3(t) 
        + C^3\ep(t). 
    \end{aligned}
    $$
    By the norm equivalence (\ref{equ-normequi}), our proof is completed.
\end{proof}

 \subsection{Estimates for density}\label{sec-ee}
 In this section, we take advantage of the wave structure of (\ref{eq-compressible}) and estimate the $L^2$-time integrability of $||\geta||_{m-1}$.
 For brevity, we rewrite (\ref{fene}) as follows
\begin{equation}\label{eq-compressible-potential}
    \begin{cases}
        \p_t \eta + \div\qu = \ret,\\
        \p_t \qu +\geta -\nu\Delta\qu =\rqut,
    \end{cases}
\end{equation}
where
$$
\ret:=-u\cdot\g \eta-\eta\div u,
$$
$$
\rqut:=  \q(\frac{1}{\rho}\div\tau )-\q ( u\cdot\g u  )- \mu\q( \er\Delta u )
-(\mu+\lambda)\q ( \er\nabla\div u ) -\q (\eta\geta).
$$
\begin{prop} \label{prop-ee}
    Suppose that the assumptions in Proposition \ref{prop-key lemma} hold, then for $t>0$, the following inequality holds:
    $$
    \begin{aligned}
         \ee(t)\leq \eb(0) 
        +C\eb(t) +  C\ep(t)   
        +C\eb^{\frac{3}{2}}(t)
        +C\ep^{\frac{3}{2}}(t)
        +C\ee^{\frac{3}{2}}(t).
    \end{aligned}
    $$
\end{prop}
\begin{proof}
    Since $\eta$ satisfies the \poin inequality, we only need to consider the highest order estimate.
    By integrating (\ref{eq-compressible-potential}) with $(\ett \gmm \g \eta,\ett \gmm\qu)$ and then integrating the resulting equality in time, we obtain
    \begin{equation}\label{eq-ee-hm}
            \it \ets||\gmm\g \eta||^2\d s=
             \it\ets||\gmm\g\qu||^2\d s+ R_1+R_2+R_3,
    \end{equation}    
    where
    $$
    \begin{aligned}
        R_1:= &
        \nu\it\ix\ets\gmm\Delta\qu\cdot\gmm\geta\dx \d s\\
        &-\it\ix\ets\p_t(\gmm\qu\cdot\gmm\geta)\dx\d s,  \\
        R_2:=& \it \ix\ets\gmm\g\ret\cdot\gmm\qu\dx\d s ,\\
        R_3:= & \it \ix\ets\gmm\rqut\cdot\gmm\geta\dx\d s.
    \end{aligned}
    $$

    For linear terms $R_1$, we use \poin inequality to obtain
    $$
    \begin{aligned}
        R_1=&
    -\ix\ett( \gmm\qu\cdot\gmm\geta )(t)\dx\\
    &+\ix( \gmm\qu\cdot\gmm\geta )(0)\dx\\
    &    + \frac{2\tc}{\nu}\it\ix\ets\gmm\qu\cdot\gmm\geta\dx\d s\\
    &+\nu \it\ix\ets\gmm\Delta\qu\cdot\gmm\geta\dx \d s\\
    \leq & \frac{1}{2}\ett||\gmm\qu||^2 
    +\frac{1}{2}\ett||\gmm\g\eta||^2\\
     &+\frac{1}{2} ||\gmm\qu(0)||^2+ \frac{1}{2}||\gm\eta(0)||^2\\
    &+\ne\it\ets||\gmp\qu||^2\d s 
    + \frac{1}{4}\it\ets||\gmm \geta||^2\d s.
    \end{aligned}
    $$

    For nonlinear terms $R_2$, we deduce from \holder inequality, Sobolev inequality and Lemma \ref{le-tau} that
    $$
    \begin{aligned}
        R_2\le & \it \ets ||u||_{m-1}||\geta||_{m-1} ||\gmm\qu||\d s
        + \it \ets ||\eta||_{m-1} ||\g\qu||_{m-1}^2\d s\\
        \le & \supt\bl ||\eta(s)||_{m-1}+ ||u(s)||_{m-1} \brr
        \it\ets\bbl ||\geta||_{m-1}^2+ ||\g\qu||_m^2\bbr\d s,
    \end{aligned}
    $$

    Using Sobolev inequality, (\ref{low bound}) and Lemma \ref{le-tau}, we can bound $R_3$ by 
    $$
    \begin{aligned}
        R_3\leq & C \it \ets ||\gm\psi||_\hyier \d s 
        + \frac{1}{4} \it \ets ||\gmm\geta||^2 \d s\\
        & +  \it \ets ||u||_{m-1}||\g u||_{m-1}||\gm\eta||\d s
         + \mu\it \ets  ||\eta||_{m-1}||\g u||_m\,||\gm\eta|| \d s\\
        & +\nu \it \ets ||\eta||_{m-1}||\g\qu||_m\,||\gm\eta|| \d s
         + \it \ets ||\eta||_{m-1}||\geta||_{m-1}^2 \d s\\
         \leq &  C \it \ets ||\gm\psi||_\hyier \d s 
        + \frac{1}{4} \it \ets ||\gmm\g \eta||^2 \d s\\
        & + C \supt \bl ||\eta(s)||_{m-1}+ ||u(s)||_{m-1}\brr
        \it \ets\bbl ||\geta||_{m-1}^2 + ||\g u||_m^2
        + \ne||\g \qu||_m^2\bbr \d s.
    \end{aligned}
    $$

    Combining the (\ref{eq-ee-hm}) with the above bounds of $R_1$, $R_2$ and $R_3$ yields
    \begin{equation*}\label{est-ee-hm}
        \begin{aligned}
            &\it \ets||\gm\eta||^2\d s\\
            \leq & C||\gmm\qu(0)||^2+ C||\gm\eta(0)||^2 +  C\ett||\gmm\qu||^2 
    +C\ett||\gm\eta||^2\\
    & + C\it \ets||\psi||_{m,\hyi}^2 \d s
    +C\ne\it\ets||\gmp\qu||^2\d s\\
    & + C \supt \bl ||\eta(s)||_{m-1}+ ||u(s)||_{m-1}\brr
        \it \ets\bbl ||\geta||_{m-1}^2 + ||\g u||_m^2
        + \ne||\g \qu||_m^2\bbr \d s.
        \end{aligned}
    \end{equation*}

Combining the above inequality with Lemma \ref {le-poincare}, we deduce that
$$
\it ||\geta||_{m-1}^2\d s\leq C\eb(0) 
        +C\eb(t) +  C\ep(t)   
        +C\eb^{\frac{3}{2}}(t)
        +C\ep^{\frac{3}{2}}(t)
        +C\ee^{\frac{3}{2}}(t).
$$
By the definition of $\ee$, our proof is completed.
\end{proof}

\section{Incompressible limit}\label{sec-incompressible limit}

    In this section, we focus on establishing the convergence rate of the incompressible limit for system (\ref{fene}). Our result indicates that, when $\lambda\rightarrow \infty$, (\ref{fene}) will converge to (\ref{equ-limit}), with the convergence rate decreasing over time. 

We first prepare the well-posedness property of (\ref{equ-limit}).
\begin{prop}\label{prop-limit system-decay}
    Suppose that $||v_0||_m+ ||\varphi_0||_{ m,\ler }\leq \varepsilon$ is the initial data of (\ref{equ-limit}). Then (\ref{equ-limit}) has a unique solution $(v,\Pi,\varphi)$ satisfying 
    $$
    ||v(t)||_m+ ||\varphi(t)||_{m,\ler}\leq C \varepsilon,
    $$
    where $C>0$ is a constant.
\end{prop}
 Using the above result,  along with the decay property in Theorem \ref{thm-exist}, we are able to prove Theorem \ref{thm-limit}.
 Note that by the proof of Theorem \ref{thm-exist}, assumption (\ref{eq-assumption}) holds as long as the solution exists. For convenience, we assume $ (\eta, u,\psi) $ and $ (v,\varphi)$ satisfy the following bounds:
\begin{equation}\label{smallness-limit part}
     \cal E(\eta(t), u(t),\psi(t)) \leq \frac{1}{64\bc^2},\quad
 ||v(t)||_m^2+ ||\varphi(t)||_{m,\ler}^2\leq\frac{1}{64 \bc^2}.
\end{equation}
\begin{proof}[Proof of Theorem \ref{thm-limit}]
    By computing $(\ref{incompressible})$-$(\ref{equ-limit})$ and denoting $U:=\pru-v$, $\Psi=\psi-\varphi$, we deduce
    \begin{equation}\label{equ-converge}
     \begin{cases}
         \p_t U - \mu \Delta U
         = - (\P \div(\rho u \otimes u) -v\cdot \g v) + \P \div \tau - \div \t\tau - \mu\Delta \P\eu,\\
         \begin{aligned}
        \p_t \Psi-\h L\Psi =& -(u\cdot\g\psi-v\cdot\g\varphi)
         -\divr ( -\g u\cdot R \psi + \g v\cdot R \varphi) \\
         &+ \divr ( -\g ( u-v)\cdot R\pin). 
         \end{aligned}
     \end{cases}
    \end{equation}
    We begin with the $L^2$-estimate. To characterize the decay rate, we set
        \begin{equation*}
        b(t)=
        \begin{cases}
            \nu  (1+ \delta t ),\qquad &0<t\leq  \nu,\\
            C_\delta \nu^2\ett,\qquad &t> \nu.
        \end{cases}
    \end{equation*}
    By integrating the above equation with $(b(t)U,b(t)\Psi)$ and using (\ref{linearpart-equal}), we obtain
    \begin{equation}
        \begin{aligned}
            &\frac{1}{2}\dt \bbl  b(t) \big( ||U||^2 
            + ||\Psi||_\lerer \big) \bbr 
            +  \mu  b(t) ||\g U||^2 
            +  b(t) ||\Psi||_\hyier
            = \sum\limits_{i=1}^5 S_i.
        \end{aligned}
    \end{equation}
    where
    $$
    \begin{aligned}
        S_1:=&\frac{1}{2}b'(t)\bbl ||U||^2 + ||\Psi||_\lerer\bbr 
            + b(t) \ixr \divr ( -\g\qu\cdot R \pin)\Psi\drx,\\
        S_2:=&\mu b(t) \ix \g\P \eu \cdot \g U\dx 
            + b(t)\ixr \divr (\P\g\eu\cdot R\pin )\Psi\drx,\\
        S_3:=& b(t)\ixr \divr (-\g u\cdot R \psi+ \g v\cdot R \varphi)\Psi\drx,\\
        S_4:=&- b(t)\ix (\P\div(\rho u \otimes u)-v\cdot\g v)\cdot U\dx,\\
        S_5:=& - b(t)\ixr (u\cdot\g \psi-v \cdot \g \varphi)\Psi\drx.
    \end{aligned}
    $$

    By (\ref{initial-additional}) and \poin inequality, we have 
    \begin{equation}\label{poin-U}
    \begin{aligned}
    ||U||^2 \leq &
    \ix\bbl\bl  U- \ix U \d x\brr
    +\ix \P\eu \d x\bbr
    ^2 \dx \\
    \leq &  2 ||\g U||^2 + C ||\eta u||^2.
    \end{aligned}
    \end{equation}
    By the choice of $\delta$ and $\nu$ in (\ref{delta}), we can bound the linear term $S_1$ by 
    \begin{equation*}\label{s1s2}
    \begin{aligned}
        S_1 \leq & \frac{1}{2}\mu  b(t) ||\g U||^2 + \frac{1}{2} b(t) ||\Psi||_\hyier \\
        &+ C b(t)\bbl 
        ||\eta||_2^2||\g\pru||^2 + ||\eta||_2^4||\g u||^2 + ||\g \qu ||^2 \bbr.
    \end{aligned}
    \end{equation*}

    For $S_2$ and $S_3$, we take advantage of (\ref{eu}) and (\ref{smallness-limit part}) to obtain
    \begin{equation*}\label{s3s4}
        \begin{aligned}
            S_2 \leq & \mu b(t) ||\g\eu||\,||\g U||
            + C b(t) ||\g\eu ||\,||\Psi||_\hyi\\
            \leq & C b(t)\bbl||\eta||_2^2 ||\g\pru||^2 + ||\eta||_2^4||\g u||^2 + ||\eta||_2^2||\g\qu||^2\bbr \\
            &+ \frac{1}{10}\mu b(t) ||\g U||^2 
            + \frac{1}{8} b(t) ||\Psi||_\hyier,\\
            S_3= & b(t) \ixr (\g U\cdot R\varphi) \cdot\gr\frac{\Psi}{\pin}\dr\dx\\ 
            &+ b(t) \ixr (\g v\cdot R \Psi) \cdot\gr\frac{\Psi}{\pin}\dr\dx\\
            & + b(t) \ixr (\g(\qu-\eu)R\psi )\cdot
            \gr\frac{\Psi}{\pin}\dr\dx\\
            \leq & C b(t) ||\g\qu||^2 
            + C b(t) ||\eta||_2^4 ||\g u||^2
            + \frac{1}{10}\mu b(t) ||\g U||^2
            + \frac{1}{4}  b(t) ||\Psi||_\hyier.
        \end{aligned}
    \end{equation*}

    To estimate $S_4$, we use (\ref{u-decompose}) and decompose $\div(\rho u\otimes u) $ as follows:
    \begin{equation*}
        \begin{aligned}
    \div(\rho u\otimes u) =&\pru\cdot\g \pru 
    + \div ((\qu+\q\eu)\otimes u)\\
    & + \div(\pru\otimes (\qu - \P \eu).
        \end{aligned}
    \end{equation*}
    Hence,
    \begin{equation*}
        \begin{aligned}
    &\div(\rho u \otimes u)-v\cdot\g v \\
    = & U\cdot\g\pru + v\cdot\g U 
    + \div ((\qu+\q\eu)\otimes u)\\
     &+ \div(\pru\otimes (\qu - \P \eu).
        \end{aligned}
    \end{equation*}
    By combining the above equation with Lemma \ref{le-poincare}, (\ref{eu}), (\ref{smallness-limit part}), (\ref{poin-U}) and $\div v=0$, one has
    \begin{equation*}\label{s6}
        \begin{aligned}
            S_4= & b(t) \ix (U\cdot\g U) \cdot \pru \dx\\
            &- b(t) \ix \div ((\qu+\q\eu)\otimes u)\cdot U \dx\\
            & - b(t) \ix \div(\pru\otimes (\qu - \P \eu)\cdot U \dx\\
            \leq & \bc b(t) ||\pru||_1(||\g U|| + C||\eta u||)||\g U||\\
            &+ C b(t) (||u||_1 + ||\pru||_1 )(||\eta u||_1 + \qu||_1)||\g U||\\
            \leq & C  b(t) \bbl ||\eta||_2^2||\g\pru||^2 + ||\eta||_2^4||\g u||^2 + ||\g\qu||^2 \bbr \\
            & + \frac{1}{5}\mu b(t) ||\g U||^2.
        \end{aligned}
    \end{equation*}
    Finally, similar as the above arguments, we can bound $S_5$ by
    \begin{equation*}\label{s5}
    \begin{aligned}
        S_5 = &-b(t)\ixr( U\cdot\g \psi) \Psi \drx \\
        &+b(t)\ixr (\P\eu\cdot\g\psi)\Psi \drx\\
        &-b(t)\ixr (\qu\cdot\g\psi)\Psi \drx\\
        \leq &
        C  b(t) \bbl ||\eta||_2^2||\g\pru||^2 + ||\eta||_2^4||\g u||^2 + ||\g\qu||^2 \bbr ||\psi||_{2,\hyi}^2\\
        &+\frac{1}{10} \mu  b(t) ||\g U||^2 
        + \frac{1}{8} b(t) ||\Psi||_\hyier.
    \end{aligned}
    \end{equation*}

   Collecting the bounds of $S_1$-$S_5$, we deduce
   \begin{equation}\label{est-UPsi}
   \begin{aligned}
   & \dt \bbl b(t) \big( ||U||^2 
     + ||\Psi||_\lerer \big) \bbr\\
    \leq & C b(t)\bbl 
    ||\eta||_2^2||\g\pru||^2 + ||\eta||_2^4||\g u||^2 + ||\g \qu ||^2 \bbr.
   \end{aligned}
   \end{equation}
   Noticing that
   $$
   \begin{aligned}
   &\it b(t) ||\eta(s)||_2^2||\g\pru(s)||^2\d s\\
   \leq & C \bbl \supt \nu\ets ||\eta(s)||_2^2\bbr 
   \it a(s) ||\pru(s)||^2 \d s\leq C \varepsilon^4,\\
    &\it  b(t) ||\eta(s)||_2^4||\g u(s)||^2\d s \\
    \leq & C \bbl \supt \ne\ets ||\eta(s)||_2^4\bbr
    \it ||\g u (s)||^2\d s\leq C \varepsilon^4 ,\\
    &\it b(t) ||\g \qu(s)||^2\d s
    \leq C \varepsilon^2, 
   \end{aligned}
   $$
   we integrate (\ref{est-UPsi}) in time and take advantage of the above estimates and (\ref{initial-additional}) to obtain
   \begin{equation}\label{est-U-l2}
          b(t) ||U(t)||^2 
     + b(t)||\Psi(t)||_\lerer 
     \leq C\varepsilon^2,
   \end{equation}
   
   Next, we consider the $\dot H^{m-1}$-estimate. By applying $\gmm$ to (\ref{equ-converge}), we deduce that $(U\mm,\Psi\mm)$ satisfies:
   \begin{equation*}
     \begin{cases}
     \begin{aligned}
         \p_t U\mm - \mu \Delta U\mm
         =& - \gmm(U\cdot\g\pru + v\cdot\g U) \\
         &+\gmm R_U
         + \P \div \tau\mm - \div \t\tau\mm,
     \end{aligned}
         \\
         \begin{aligned}
         \p_t \Psi\mm-\h L\Psi\mm = 
         &- \gmm (U\cdot\g\psi + v\cdot \g \Psi)
         -  \gmm \divr( \g U \cdot R \psi
         + \g v\cdot R \Psi )\\
         &+ \divr ( -\g U\mm \cdot R\pin)
         + \gmm R_\Psi. 
         \end{aligned}
     \end{cases}       
   \end{equation*}
   where
   $$
   \begin{aligned}
    R_U:=&
    - \div ((\qu+\q\eu)\otimes u)
     - \div(\pru\otimes (\qu - \P \eu)
      - \mu\Delta \P\eu\\
   R_\Psi :=& \P\eu \cdot\g\psi 
   - \qu\cdot\g \psi
   + \divr (\g\P\eu \cdot R \psi)
   + \divr ( -\g\qu\cdot R\psi)\\
   & + \divr ( \g\P\eu \cdot R \pin ) 
   + \divr ( - \g\qu \cdot R \pin ). 
   \end{aligned}
   $$

   By dotting the above system with $(b(t) U\mm, b(t) \Psi\mm)$, one has
   \begin{equation}\label{equ-UPsi-hmm}
       \begin{aligned}
           \frac{1}{2}\dt 
           \bbl b(t) \big( ||U\mm||^2 +  ||\Psi\mm||_\lerer \big)  \bbr
           + \mu b(t)||\g U \mm||^2 
           + b(t) ||\Psi\mm||_\hyier
           =\sum_{i=1}^6 T_i.
       \end{aligned}
   \end{equation}
   where
   $$
   \begin{aligned}
       T_1:= & \frac{1}{2}b'(t)
           \bbl ||U\mm||^2 + b(t) ||\Psi\mm||_\lerer \bbr, \\
       T_2:= & - b(t) \ix \gmm(U\cdot\g\pru + v\cdot\g U) \cdot U\mm \dx,\\ 
       T_3:=&- b(t) \ixr \gmm (U\cdot\g\psi + v\cdot \g \Psi)\Psi\mm \drx\\
       &-  b(t) \ixr \gmm \divr( \g U \cdot R \psi
         + \g v\cdot R \Psi ) \Psi \mm \drx, \\
       T_4:= &  b(t) \ix \gmm R_U\cdot U\mm dx,\\
       T_5:=&  b(t) \ixr \gmm R_\Psi \Psi \mm \drx.
   \end{aligned}
   $$

   For the linear term $T_1$, we use (\ref{delta}), \poin inequality and Lemma \ref{le-poincare-r} to obtain
   \begin{equation*}\label{t1}
       T_1\leq \frac{1}{4}b(t)
       \bbl \mu ||\g U\mm||^2 + ||\Psi\mm||_\hyier\bbr.
   \end{equation*}

   For $T_2$, by $\div v=0$, Lemma \ref{le-communicator}, (\ref{eu}), (\ref{smallness-limit part}),  (\ref{poin-U}) and \poin inequality, one has
   \begin{equation*}\label{t2}
   \begin{aligned}
       T_2 = & b(t) \ix \g^{m-2} (U\cdot \g \pru)\cdot U^m \dx \\
       &-b(t) \ix ([\gmm,v\cdot\g] U )\cdot U\mm \dx\\
       \leq &b(t) ( \bc ||\g U \mm||+ C||\eta u|| )
       ||\g\pru||_{m-1} ||\g U \mm||\\
      & + b(t) \bc ||v||_{m}||\g U\mm||^2\\
       \leq & C b(t) ||\eta||_1^2||u||_1^2||\g\pru||_{m-1}^2 
        +   \frac{1}{2} \mu b(t) ||\g U\mm||^2.
   \end{aligned}
   \end{equation*}

   Similarly, $T_3$ can be bounded by:
       \begin{align*}
           T_3 \leq &
           C b(t) ||\eta||_1^2||u||_1^2||\psi||_{m,\hyi}^2 
           + \frac{1}{8}\mu b(t) ||\g U\mm||^2 
           + \frac{1}{2} b(t) ||\Psi\mm||_\hyier.
       \end{align*}

       It remains to estimate $T_4$ and $T_5$. By Lemma \ref{le-poincare}, \ref{le-poincare-r} and (\ref{eu-hs}), (\ref{smallness-limit part}), we obtain
       \begin{equation*}
           \begin{aligned}\label{t5}
               T_4\leq & C b(t)||(\qu,\eu)||_{m-1}||u||_{m-1}||\g U\mm||\\
               & + C b(t) ||\pru||_{m-1}||(\qu,\eu)||_{m-1}||\g U\mm||\\
               &+ C b(t) ||\eta u||_{m}||\g U \mm||\\
               \leq & C b(t) \bbl
               ||\eta||_m^2 ||\g\pru||_{m-1}^2
               + ||\eta||_m^4 ||\g u||_{m-1}^2
               + ||\g\qu||_{m-1}^2 \bbr\\
               &+ \frac{1}{8}\mu b(t) ||\g U\mm||^2,\\
               T_5\leq & C b(t) ||\eta u||_{m-1}||\psi||_{m.\hyi}||\Psi\mm||_\hyi\\
               &+ C b(t) ||\psi||_{m,\ler} ||\g\qu||_{m-1}||\Psi\mm||_\hyi\\
               & + C b(t) (||\psi||_{m-1,\ler}+1) ||(\g\eu,\g\qu)||_{m-1} ||\Psi\mm||_\hyi\\
               \leq & C b(t) \bbl
               ||\eta||_m^2 ||\g\pru||_{m-1}^2
               + ||\eta||_m^4 ||\g u||_{m-1}^2
               + ||\g\qu||_{m-1}^2 \bbr\\
               &+ \frac{1}{4} b(t) ||\Psi\mm||_\hyi^2.
           \end{aligned}
       \end{equation*}

       Inserting the bounds of $T_1$-$T_5$ into (\ref{equ-UPsi-hmm}) and integrating the resulting inequality in time, we find 
       $$
       \begin{aligned}
       &b(t) ||U\mm(t)||^2 
     + b(t)||\Psi\mm(t)||_\lerer \\
     \leq &||U\mm(0)|| + ||\Psi\mm(0)||_\lerer\\
     &+ C \it b(s) \bbl
               ||\eta||_m^2 ||\g\pru||_{m-1}^2
               + ||\eta||_m^4 ||\g u||_{m-1}^2
               + ||\g\qu||_{m-1}^2 \bbr \d s.
       \end{aligned}
       $$
       By (\ref{initial-additional}) and the argument parallel to the corresponding part in $L^2$-estimate, we deduce
       $$
       b(t) ||U\mm(t)||^2 
     + b(t)||\Psi\mm(t)||_\lerer \leq C\varepsilon^2.
       $$
       We can obtain the desired result by combining the above inequality with (\ref{est-U-l2}).
    \end{proof}

\section*{Acknowledgment}
Jincheng Gao was partially supported by the National Key Research and Development Program of China(2021YFA1002100), Guangdong Special Support Project (2023TQ07A961), Guangzhou Science and Technology Program(2024A04J6410) and  Guangdong Basic and Applied Basic Research Foundation(2021B1515310003).
Jiahong Wu was partially supported by the National Science Foundation of the United States under DMS 2104682 and DMS 2309748. 
Zheng-an Yao was partially supported by
Shenzhen Science and Technology Program( CJGJZD2021040809140300), 
National Key Research and Development Program of China(2020YFA0712500) and R\&D project of Pazhou Lab (Huangpu)( Grant2023K0601).

\section*{Data availability}
Data sharing is not applicable to this article as no data sets were generated or analyzed during
the current study.

\section*{Declarations}
\subsubsection*{Conflict of interest}
The authors declare that they have no Conflict of interest.

\bibliography{FENE_limit}
\end{document}